\documentclass[a4paper,11pt]{article}
\author{Alan Hammond, Yuval Peres}
\title{Fluctuation of planar Brownian loop capturing large area}

\usepackage{amsfonts}
\usepackage[dvips]{graphicx}

\def\vol{{\rm vol}}
\def\enc{{\rm enc}}
\def\ext{{\rm ext}}
\def\conv{{\rm conv}}
\def\arcl{{\rm arcl}}
\def\rad{{\rm rad}}
\def\rin{{\rm R}_{\rm in}}
\def\rout{{\rm R}_{\rm out}}
\def\intn{{\rm int}}
\def\disk{{disk}}
\newtheorem{definition}{Definition}
\newtheorem{lemma}{Lemma}
\newtheorem{corollary}{Corollary}
\newtheorem{theorem}{Theorem}
\newtheorem{prop}{Proposition}

\begin{document}
\maketitle
\noindent{\bf Abstract.} We consider a planar
Brownian loop $B$ that is run for a time $T$ and conditioned on the
event that its
range encloses the unusually high area of $\pi T^2$, with $T \in
(0,\infty)$ being large.
The conditioned process, denoted by $X$, was proposed by Senya
Shlosman as 
a model for the
fluctuation of a phase boundary. We study the deviation of the range
of $X$ from a circle of radius $T$.
This deviation is measured by the inradius $\rin(X)$ and
outradius $\rout(X)$, which are the maximal radius of a {\disk} enclosed by
the range of $X$, and the minimal radius of a {\disk} that contains this range.
We prove that, in a typical realization of the conditioned measure,
each of these quantities differs from $T$ by at most $T^{2/3 + \epsilon}$.  
\begin{section}{Introduction}
The goal of this paper is to analyse the fluctuations of a
planar Brownian loop under the condition that it
encircles a large area. Throughout, $B: [0,T] \to \mathbb{R}^2$
will denote a standard planar Brownian loop, that is, a planar
Brownian motion with initial location $B(0)=0$ that is conditioned on
the event that $B(T)=0$. Allowing $\enc(B)$ to
denote the random set of points that lie in the union of all bounded 
components of $\mathbb{R}^2 \setminus \big\{ B(t) : t \in [0,T] \big\}$, our
conditioning takes the form  
\[
\big\vert
\enc (B) \big\vert \geq \pi T^2,
\] 
where $\vert \cdot \vert$ denotes
two-dimensional Lebesgue measure. Note that by the spatial-temporal
scaling satisfied by Brownian motion, the law of the conditioned
process is the same as that obtained from sampling a Brownian loop run for a
unit of time that is conditioned to enclose an area exceeding
$\pi T$, and then dilating space by a factor of
$\sqrt{T}$. Throughout, we will define the process on the interval
$[0,T]$.
The conditioned process will be
denoted by $X : [0,T] \to \mathbb{R}^2$. 
As we discuss in Section \ref{sectwo}, a classical variational principle
suggests that the 
range $X[0,T]$ takes a form close to that of a circle of radius $T$.
The principal aim of this paper is  
to investigate the magnitude of the deviation of the range $X[0,T]$
from such a circle.
\begin{figure}
\begin{center}
\includegraphics[width=0.5\textwidth]{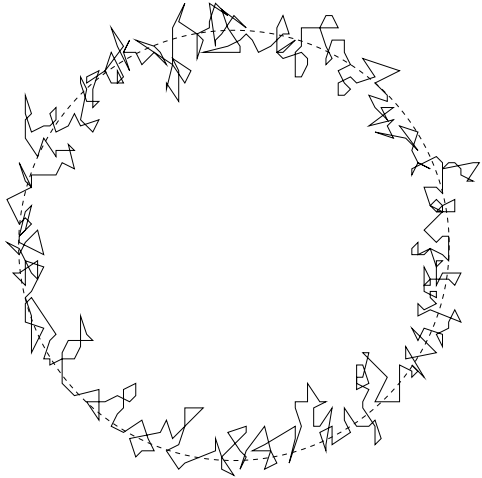} \\
\end{center}
\caption{A sketch of the conditioned motion.}
\end{figure}
Our main theorem provides a bound on a quantity that measures this
deviation. 
To be precise, for a planar compact set $K$, the
inradius  $\rin(K)$ of $K$ is the maximal radius of a circle lying
in $K$, while the outradius $\rout(K)$ is the minimal radius of any
circle in which $K$ is contained. We will write $\rin(B)$ and
$\rout(B)$ for $\rin\big( \enc(B) \big)$ and $\rout\big( \enc(B) \big)$, adopting the
same shorthand for the process $X$.
\begin{theorem}\label{tgg}
 Let $\epsilon \in (0,1/6)$. For  $c \in \big( 0,\pi^2/{32}
\big)$ 
 and all $T \geq T_c$ sufficiently high,
\begin{displaymath}
 \mathbb{P} \Big( \rin \big( X \big) < T - T^{\frac{2}{3}+\epsilon}
 \Big) \leq \exp\Big\{ - c T^{\frac{1}{3} + 2 \epsilon} \Big\},
\end{displaymath}
For any constant $\hat{c}$  satisfying $\hat{c} \in \big( 0,\pi^2/{2^9}
\big)$ and  all $T \geq T_{\hat{c}}$ sufficiently high,
\begin{displaymath}
 \mathbb{P} \Big( \rout \big( X \big) > T + T^{\frac{2}{3}+\epsilon}
 \Big) \leq \exp\Big\{ - \hat{c} T^{\frac{1}{3} + 2 \epsilon} \Big\}.
\end{displaymath}
\end{theorem}
How close is $X[0,T]$ to the boundary of its convex hull? This is a
question about the local nature of the deviation of the conditioned
process. We will write $\mathcal{L}(\conv K)$ for the length of the 
longest line segment that lies in the boundary of the convex hull
$\conv(K)$ of $K$. We also define the maximum local roughness {\rm MLR}(K) of
$K$ to be the maximal
distance between a point in $K$, and the boundary of $\conv(K)$. That
is,
\begin{displaymath}
 {\rm MLR} \big( K \big) := \sup_{k \in K} \inf_{x \in \partial (\conv(K))} d(x,k).
\end{displaymath}
We will write $\mathcal{L}(\conv B)$ for $\mathcal{L}(\conv K)$ in the case where $K$ is the range
of the process $B: [0,T] \to \mathbb{R}^2$. A similar convention will
apply for the maximum local roughness, and for the conditioned process $X: [0,T] \to \mathbb{R}^2$. 
 
Senya Shlosman proposed this model to us, presenting some heuristic
arguments that its deviation
behaviour has much in common with that observed in numerous models of
phase boundaries in two-dimensional random systems, more specifically,
that exponents describing the typical behaviour of $\mathcal{L}(\conv X)$ and
${\rm MLR}(X)$ coincide with those in these other models. We will 
present some 
heuristic arguments of our own in favour of this belief at the end of the introduction.
Theorem \ref{tgg} has the following straightforward consequence.
\begin{corollary}\label{nthmone}
 The fluctuation of the conditioned process $X : [0,T] \to
 \mathbb{R}^2$ satisfies the following bounds. 
For any $\epsilon \in (0,1/{12})$,  $c \in \big(0, 2^{-21} \big)$ 
and for all $T \geq T_c$ sufficiently high,
\begin{displaymath}
  \mathbb{P} \Big( \mathcal{L} \big( X \big) \geq T^{\frac{5}{6} + \epsilon}
  \Big) \leq \exp \Big\{ - c T^{\frac{1}{3} + 4 \epsilon } \Big\},
\end{displaymath} 
Moreover, for $\epsilon \in (0,1/6)$ and all $T \geq T_{\hat{c}}$ sufficiently high,
\begin{displaymath}
 \mathbb{P} \Big( \arcl \big( \partial ( \conv X ) \big) > 2\pi \big( T
 + T^{\frac{2}{3} + \epsilon} \big) \Big) \leq  \exp\Big\{ - \hat{c} T^{\frac{1}{3} + 2 \epsilon} \Big\}, 
\end{displaymath}
where the constant $\hat{c}$ appears in Theorem \ref{tgg}, and where
$\arcl$ denotes the arclength of a planar set.
\end{corollary}
\noindent{{\bf An outline of the proof.}}
We now discuss the techniques required to prove Theorem \ref{tgg},
outlining the structure of the paper as we do so. The first part of Section \ref{sectwo}
describes more precisely the assertion from the
theory of large deviations that the range $X[0,T]$ resembles a circle
of radius $T$. The second part is devoted to developing the tools required to
prove Theorem \ref{tgg}, while the proofs are given
in Section \ref{secthr}.  

We seek to understand the behaviour of the conditioned
motion $X$ by considering the polygon $P$ whose vertices are the locations
of $X$ at $m$ equally spaced moments of time. We consider the
area of this polygon and that of its convex hull, as well as the area
trapped between the range of the motion as it traverses the space
between two successive vertices of the polygon, and the line segment
between this pair of vertices. We are free to choose the value of $m$
as we please, and it is of little surprise given our belief about the
true fluctuation in the model that choices close to $T^{1/3}$ are
convenient.

The two most significant results in the second part of Section \ref{sectwo} are  
Lemma \ref{lemcl} and Proposition \ref{ttwo}. 
The former shows that the convex hull of the polygon $P$ is likely to
trap an area that is not much less than that captured by the
conditioned motion, for values of $m$ slightly less than
$T^{1/3}$. The discrepancy is shown to be at most a little more than
$T^{4/3}$ with high probability. Proposition \ref{ttwo} provides an estimate on the regularity of the
conditioned motion, to the effect that it is unlikely to move too quickly in short
periods of time. More precisely, in a time of order $T^{2/3 +
\epsilon}$, we rarely see the motion cover as much distance as $T^{2/3
+ 2 \epsilon}$. 
To prove each of these results, the first step is to provide a lower
bound on the probability that the Brownian loop $B : [0,T] \to
\mathbb{R}^2$ in fact satisfies the requirement of the conditioning,
that it captures an area of $\pi T^2$. This bound is provided in Lemma 
 \ref{ltwo} by estimating the probability that a regular polygon with
 an order of $T^{1/3}$ vertices and having this area is enclosed by
 the motion $B$.
We then prove the two results by showing that if $\vert \conv(P)
\vert$ is less than $\vert \enc(X) \vert$ by $T^{4/3 + \epsilon}$, or if $X$ does
move a distance of $T^{2/3+2\epsilon}$ in some interval of time of
order $T^{2/3 + \epsilon}$, then certain functionals  
 of collections of normally distributed random
variables assume high values. For example, in the proof of Lemma
\ref{lemcl}, we consider the event that $\vert \conv(P)
\big \vert$ is less than  $\vert \enc(X) \vert$ by at least $T^{4/3 + \epsilon}$. In Lemma \ref{assert}, we show
that the discrepancy $\enc(X) \setminus \conv(P)$ is contained in the
union over the polygonal edges $l$ of regions that are, roughly speaking,
rectangles whose long axis is $l$ and whose width is  the orthogonal
fluctuation of the motion $X$ during
the interval of time in which it traverses the edge $l$.
 
These edge lengths and fluctuations are two
collections of random variables whose distributions are readily
bounded above by some having a normal distribution. We have seen that
the sum of the products of edge lengths and corresponding fluctuations
is at least $T^{4/3 + \epsilon}$ if the event in question occurs. The square
of this
expression may be bounded above by the product of two random variables having the
$\chi^2$-distribution, by use of the Cauchy-Schwarz
inequality. Computations that make use of the formula for the density of
this distribution show that the event in question is less probable
than that of the conditioning being satisfied, by a comparison with the lower
bound given in Lemma \ref{ltwo}.  

In the proof of Theorem \ref{tgg}, we work with the polygon $P$, with
$m$ being chosen to be of order $T^{1/3}$. The convex hull of $P$ is
known typically to trap a high area by Lemma \ref{lemcl}. An
isoperimetric result, Lemma \ref{lone}, gives a lower bound on the
arclength of a planar convex body in terms of its area and its global
deviation, the latter giving rise to an excess over that occuring in
the extremal case of a {\disk}. We deduce that the global deviation of
the boundary of the convex hull of $P$ is not too high: for otherwise,
the sum of the edge-lengths of the polygon being high implies that the
$\chi^2$-distributed sum of their squares is improbably large. This
arclength is also forced to be high if there is any vertex $p$ of $P$ that
is too distant from the line segment in $\conv(P)$ which the motion is
traversing when it reaches $p$. The fluctuation of the motion between
two successive vertices of the polygon has been bounded in Proposition
\ref{ttwo}. We have obtained enough control on the motion to deduce
that its global deviation typically satisfies an upper bound whose
order is comparable to that satisfied by the deviation of the
approximating polygon, that is, little more than $T^{2/3}$.

We so demonstrate that the motion is likely to be trapped between circles
whose radii differ by an order that slightly exceeds that of
$T^{2/3}$. The bounds on ${\rm MLR} \big( X \big)$ and $\mathcal{L} \big( X \big)$
given in Corollary \ref{nthmone} are straightforward consequences. \\
{\bf Remark} Throughout, 
any time parameter takes a value in $[0,T]$. From time to time, terms
such as $t + T^{2/3}$ appear that may be greater than $T$. In such
cases, we are referring to the value on $[0,T]$ that is the value
stated reduced modulo $T$.
\subsection{Comparisons and heuristic arguments}
It is believed that a variety of models of phase boundaries in
two-dimensional random systems exhibit the same power-law
fluctuations, even though their macroscopic profiles differ. An
important example of such a model 
is that of a large
finite cluster in the supercritical phase of the site percolation
model in $\mathbb{Z}^2$. Choosing a parameter value $p > p_c$, 
Alexander and Uzun \cite{alexuzun} condition on the event that there
exists an open dual circuit surrounding the origin and enclosing an
area of at least $n^2$, for large $n$.
The 
asymptotic shape of this circuit is the boundary of a compact convex body, known as
the Wulff crystal, that minimises a surface tension,
c.f. \cite{wulff}.
The fluctuation of the circuit away from this shape 
may be measured by the maximum local roughness, which means in this
case, the maximum distance of a vertex in the circuit from the
boundary of its convex hull. 
In \cite{alex}, the maximum local roughness
is established to be bounded
above by a quantity of the order of $n^{2/3}$. 
The average local roughness (which is roughly speaking the mean
distance of a vertex in $C$ from the convex hull)
is bounded above by $n^{1/3}$ up to power
order, a bound which is believed to be sharp. In \cite{alexuzun}, the maximum local roughness is shown to
satisfy a lower bound that is given by $n^{1/3}$ if we omit logarthmic corrections.

We expect that the exponents describing the typical behaviour of the
two measures of fluctuation, $\mathcal{L}$ and ${\rm MLR}$, coincide with those
anticipated for the percolation problem. 
That
is, we expect that $\mathcal{L}(\conv X)$ behaves as $T^{2/3}$, and ${\rm MLR}(X)$ as
$T^{1/3}$. To give an argument that supports the claim that $\mathcal{L}(\conv X)$
is typically not much greater than $T^{2/3}$, suppose that we sample
the measure $X$ and find a realization $X(\omega)$ where there is a line segment
$L = \big[ x_1,x_2 \big]$ in $\partial \big( \conv(X) \big)$ whose
length exceeds $T^{2/3 + \epsilon}$. The times $t_1$ and $t_2$ at
which $X(\omega)$ visits the endpoints $x_1$ and $x_2$ presumably satisfy
$\big\vert t_1 - t_2 \big\vert > T^{2/3 + \epsilon/2}$, since the
conditioned motion tends to move at a fairly constant rate. (Indeed,
in Proposition \ref{ttwo}, we will prove that the process $X$ is unlikely to
cover distances as big as $L$ at speeds significantly greater than
the average one at which $X$ moves). Choosing two points $t$ and
$t^*$ on the interval $[0,T]$ uniformly and independently of other
randomness, we may resample the path of $X$ on the interval
$\big[t,t^* \big]$, replacing $X [t,t^*]$ by a Brownian bridge
that moves from $X(t)$ to $X(t^*)$ in time $t^* - t$. The
Markov chain on loops that performs this resampling and jumps to the
new path provided that it captures the required area of $\pi T^2$, and
stays put in the other case, has the law of the conditioned process
$X$ as its invariant measure. We see that, for the action of this
resampling on $X(\omega)$, if the points $t$ and $t^*$ happen to be
picked near to $t_1$ and $t_2$ respectively, then the effect of the
resampling is to replace the motion of $X(w)$ as it traverses $L$ by
a new motion. This motion has a time of at least $T^{2/3 +
\epsilon/2}$ to traverse a distance of $T^{2/3 + \epsilon}$. 
This new section of path typically fluctuates  orthogonally to $L$ by
a distance at least of order $T^{1/3 + \epsilon/4}$ (this being
the square root of the available time). 
With a
probability that is uniformly bounded below in $T$, this fluctuation
occurs for a fixed but high fraction of time in the direction away from the existing convex
hull of $X$. In this case,
the resampled motion would
seem to capture
an area of the plane that exceeds that captured by $X(w)$ by an amount
of the order of
$ T^{2/3 + \epsilon} \cdot T^{1/3 + \epsilon/4}  = T^{1 + 5 \epsilon/4}$ (the left-hand-side here being the product of
the length of $L$ and the orthogonal fluctuation of the resampled motion).   
In this event, the resampling certainly meets the area criterion. 
Figure \ref{resam} shows a sketch of the range of a typical
realization of $X$, and a resampling that creates more trapped area by
the means just described.
\begin{figure}\label{resam}
\begin{center}
\includegraphics[width=0.75\textwidth]{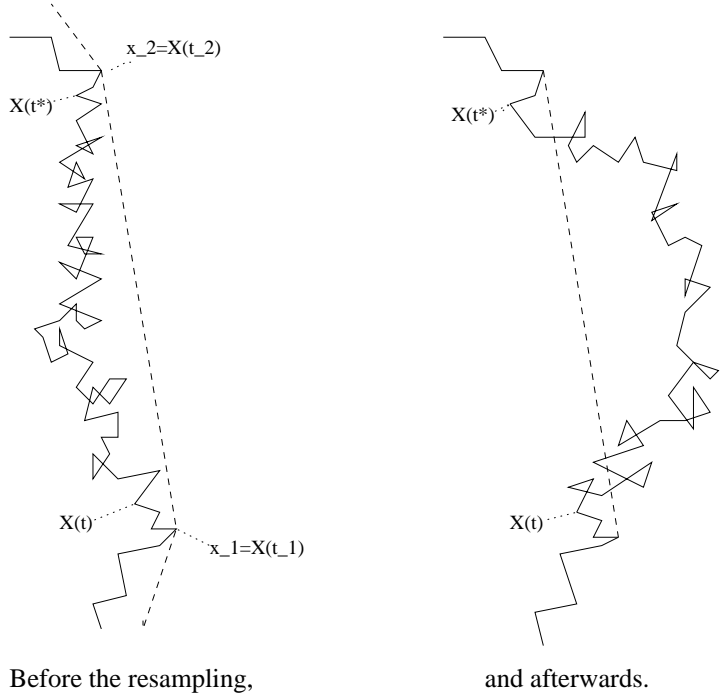} \\
\end{center}
\caption{How resampling may trap more area.}
\end{figure}
It is very
believable that the typical order of the excess of area that $X$
captures over what it must capture, $\big\vert
\enc(X) \big\vert - \pi T^2$, is linear in $T$, and that
\begin{equation}\label{phf}
\mathbb{P} \Big( \big\vert \enc(X) \big\vert - \pi T^2 \geq T^{1 + \alpha} \Big)
\end{equation}
decays at a super-polynomial rate, for any given $\alpha  > 0$.
However, the preceding argument suggests that the resampled motion -
whose law is that of $X$ - has an excess of area of order of $T^{1 +
5\epsilon/4}$ with a probability that is at least a polynomially
decaying multiple of the probability that there exists a line segment
in $\partial \big( \conv (X) \big)$ of length $T^{2/3 + \epsilon}$ (the
fact that the points $t$ and $t^*$ must be chosen to be near $t_1$
and $t_2$ is responsible for the appearance of a polynomial factor
here). So, one expects that the probability of such a line segment in a
realization of $X$ decays at a super-polynomial rate in $T$. 
It remains an interesting problem to derive such an upper bound on the quantity
(\ref{phf}), as does that of obtaining lower bounds on
the area captured after resampling (an example of the difficulties
involved in determining the area captured is
the fact that Lemma \ref{assert} is not valid if the instance of $\conv(P)$
in its statement is replaced by $P$).
\end{section}
\begin{section}{Understanding how the process $X$ fluctuates}\label{sectwo}
\begin{subsection}{The macroscopic profile of the range of
$X$}\label{sectwoone}
The theory of large deviations is of use in deriving 
the asymptotic shape of the range of the process $X$. Let
$C_0([0,1],\mathbb{R}^2)$ denote the space of continuous planar-valued
functions $f:[0,1]\to \mathbb{R}^2$ for which $f(0)=0$, and let $\it{W}$
denote two-dimensional Wiener measure. Writing $f_T(\cdot) =
T f(\cdot/T)$, we have, by Theorem 5.1
of \cite{varadhan},
\begin{equation}\label{upbdd}
 \limsup_{T \to \infty}{\frac{1}{T} \log {\it W} \Big\{ g:
 g\Big\vert_{[0,T]} = f_T, \, \, \textrm{for some $f \in C$} \Big\}}
 \leq - \inf_{x \in C}{I(x)}
\end{equation}
for $C \subseteq C_0([0,1],\mathbb{R}^2)$ closed, and
\begin{equation}\label{lowbdd}
 \liminf_{T \to \infty}{\frac{1}{T} \log {\it W} \Big\{ g:
 g\Big\vert_{[0,T]} = f_T, \, \,  \textrm{for some $f \in O$} \Big\}}
 \geq - \inf_{x \in O}{I(x)}
\end{equation}
for $O \subseteq C_0([0,1],\mathbb{R}^2)$ open.

The {\it large deviations' rate function} $I:  C_0([0,1],\mathbb{R}^2)
\to [0,\infty]$ is given by
\[
 I(f) = \frac{1}{2} \int_0^1  \vert\vert D f \vert\vert^2 (t) dt,
\]
if $f \in H_1 \big( [0,1], \mathbb{R}^2\big)$ (that is, if $f$ is absolutely continuous with square integrable derivative
$D f$), with $I(f)=\infty$ otherwise. 
By applying the contraction principle of large deviations \cite[Section
4.2.1]{Dembozeit} to the mapping of the space of Brownian paths
$Z:[0,T] \to \mathbb{R}^2$ to the space of Brownian loops $B:[0,T] \to
\mathbb{R}^2$ given by $Z(t) = B(t) - \frac{t}{T}B(T)$, we learn that
(\ref{upbdd}) and (\ref{lowbdd}) are valid for the measure $dB$
provided that the space $C_0([0,1],\mathbb{R}^2)$
is replaced by its subspace $\overline{C}_0$ consisting of functions $f$ for which
$f(1)=0$.
In evaluating the area enclosed by a loop, we will use for the present
argument the signed area, given by
\begin{equation}\label{signedarea}
A(f) = \frac{1}{2} \int_0^1 \big(
f_1(t)f_2'(t)-f_1'(t)f_2(t)\big)dt,
\end{equation}
for any $f \in H_1\big([0,1], \mathbb{R}^2 \big)$.
Noting that $A(f_T) \geq \pi T^2$ for such $f$ 
if and only if  $A(f) \geq \pi$, 
we now identify those functions in $\{ f \in
 H_1\big([0,1], \mathbb{R}^2 \big): A(f)  \geq \pi \}$ 
that minimise $I(f)$. Any $f$ in
this set is certainly square integrable, and
thus, has an $L^2$-convergent complex Fourier series
\begin{equation}\label{ffourier}
 f(x) = \sum_{n \in \mathbb{Z}}{a_n \exp \{ 2  n \pi i x \} }.
\end{equation} 
Note that, in these terms,
\begin{equation}\label{formi}
I(f) = 2 \pi^2 \sum_{n \in \mathbb{Z}}{n^2 \vert a_n\vert^2}. 
\end{equation}
The formula
 (\ref{signedarea}) for signed area translates to   
\begin{equation}\label{forma}
 A(f) = \pi \sum_{n \in \mathbb{Z}}{n \, \vert a_n \vert^2}
 = \pi \sum_{n = 1}^{\infty}{n \Big( \vert a_n \vert^2 - \vert a_{-n}
 \vert^2 \Big)}
,
\end{equation}
It is clear from 
(\ref{formi}) and (\ref{forma}) that any  $f:[0,1] \to \mathbb{C}$
that minimises $I(f)$ among those functions for which 
$A(f) \geq \pi$, and $f(0) = f(1)$,
has $a_n = 0$ for $n < 0$. It also follows from 
(\ref{formi}) and (\ref{forma}), that 
\begin{equation}\label{moref}
I(f) \geq 2\pi A(f),
\end{equation}
for $f$ having only positive Fourier modes. 
Thus, we
must have $I(f) \geq 2\pi^2$ for functions $f$ such that $A(f) \geq \pi$. 
Note however that if $f(x) = -a + a \exp \big\{ 2 \pi i x
\big\}$, for $a \in \mathbb{C}$ such that $\vert a \vert = 1$, then
$I(f) = 2 \pi^2$ and $A(f) = \pi$, so that equality in (\ref{moref})
is attained for such functions $f$. Noting that, if $f$ has only
positive Fourier modes and has some $a_n \not= 0$ for $n > 1$, the
inequality in (\ref{moref}) is strict, we have deduced that
each of the minimising functions takes the form 
of a progression at constant rate along the circumference of a circle
of radius $1$.

However, it appears that the fluctuation behaviour 
of the conditioned process may not be understood
by a direct application of the techniques of large
deviations. We now begin to develop the tools required for our study
of this deviation.
\end{subsection}
\begin{subsection}{Tools for the proof}\label{sectwotwo}
Firstly, we find a lower bound on the probability that the Brownian
loop captures the required area.
\begin{lemma}\label{ltwo}
Let $G_m = G_m(T)$ denote a regular polygon with $m$ vertices that
contains an area equal to $\pi T^2$. Then, for any constant $C_1$
satisfying $C_1 > 4/3$, and for all $T \geq T_{C_1}$ sufficiently high,
\begin{eqnarray}
 & & \mathbb{P} \Big( \enc \big( B \big) \supset x + G_{\lfloor T^{\frac{1}{3}}
\rfloor}  \ \textrm{for some $x \in \mathbb{R}^2$} \Big) \nonumber \\
 & \geq & \exp \Big\{- 2
\pi^2 T - C_1 T^{\frac{1}{3}} \log T \Big\}. \nonumber
\end{eqnarray}
\end{lemma}
{\bf Proof}
Denote the successive
vertices of $G_m(T)$ by
$$
\Big\{ y_i : i \in \{1,\ldots, m \} \Big\},
$$
where we set $m= \lfloor T^{1/3}\rfloor$. 
and let $D_i$ denote the line segment whose endpoints are the centre
of  $G_{m}$, and $y_i$. Let $\phi_i$ denote the
{\disk} of radius one whose centre lies on the continuation of $D_i$ at
distance two from $y_i$. Let $M_i$ denote the open half-plane, disjoint
from $G_{m}(T)$, whose boundary contains the
line segment $[y_i,y_{i+1}]$ (in the case where $i = m$, 
the line segment $[y_{m} , y_1]$ ).
Let $q \in
\mathbb{R}^2$ be such that $q \in - \phi_1$. For  $ i \in \{0,\ldots,
m - 2\}$, let $A_i$ denote the
event that
\begin{eqnarray}
B \Big( \frac{iT}{ m} \Big) & \in & q +
\phi_{i+1}, \nonumber \\
B \Big( \frac{(i+1)T}{ m} \Big) & \in & q +
\phi_{i+2}, \nonumber
\end{eqnarray}
and
\begin{displaymath}
 B(t)  \in  q + M_{i+1} \ \textrm{for $t \in
 \Big[\frac{iT}{m} , \frac{\big(i + 1
 \big) T}{m} \Big]$}.
\end{displaymath}
In the case where $i =  m-1$, we use the same
definition, with $\phi_1$ replacing $\phi_{i+2}$ in its statement.
\begin{figure}
\begin{center}
\includegraphics[width=0.5\textwidth]{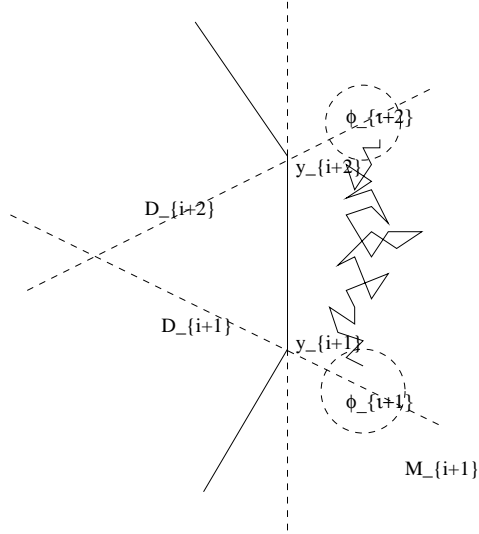} \\
\end{center}
\caption{A realization of the event $A_i$.}
\end{figure}
Note that 
\begin{equation}\label{useit}
\bigcap_{i=0}^{m-1}{A_i}  \subseteq 
\Big\{ \enc(B) \supset  q + G_{m} \Big\}. 
\end{equation}
We claim that
\begin{equation}\label{clc}
 \mathbb{P} \bigg( \bigcap_{i=1}^{ m}{
 \Big\{ B \Big( \frac{iT}{ m} \Big) \in q +
\phi_{i+1}  \Big\}} \bigg)  \geq \exp{ \Big\{- 2 \pi^2 T - \big( 2/3 +
o(1) \big)  T^{\frac{1}{3}}
\log T \Big\} }.
\end{equation}
To see this, note that the left-hand-side of (\ref{clc}) is given by  
\begin{equation}\label{llt}
 \frac{m^m}{ \big( 2 \pi T  \big)^{m - 1 } }
 \int \exp{ \Big\{ - \frac{
m}{2 T } \sum_{i=0}^{m -
1}{\vert\vert x_{i+1} - x_i \vert\vert^2}
\Big\} } d x_1 \ldots d x_{m - 1 },
\end{equation}
where the range of integration is equal to $\big( q + \phi_2 \big) \times \ldots \times 
\big( q + \phi_{m } \big)$
and where we set $x_0 = x_{m} = 0$. The form for the expression in (\ref{llt}) occurs by
computing the density of a finite-dimensional distribution of a
Brownian bridge as the ratio of the corresponding density for a
Brownian motion and the density at zero of a normal random variable
with mean zero and variance $T$. 
It is straightforward to show that 
\begin{equation}\label{pert}
 \big\vert D_i \big\vert = T + \frac{\pi^2}{3}T^{\frac{1}{3}} + O(1).
\end{equation}
Thus, the distance between successive vertices satisfies 
\begin{eqnarray}
 d \big( y_i, y_{i+1} \big) & = & 2 \big\vert D_i \big\vert \sin{
 \frac{\pi}{m} } \leq \frac{2 \big\vert D_i
 \big\vert \pi}{m} \nonumber \\
 & \leq & \frac{2 \pi T}{m} + O(1). \label{pret}
\end{eqnarray}
From (\ref{pret}), it follows that the expression in (\ref{llt}) is bounded below by 
$$
 \frac{m^m}{ \big( 2 \pi T  \big)^{m - 1 } }
 \Big[ \prod_{i = 2}^{m}{\vol \big( \phi_i \big)} \Big] \exp{ \Big\{ \frac{
 - m^2}{2T} \Big( \frac{2 \pi T}{m} + O(1) \Big)^2  \Big\}}.
$$
Since $\vol (\phi_i) = \pi$ for $ i \in \{1,\ldots, \lfloor T^{1/3}
\rfloor \}$, this expression is bounded below by 
$
 \exp{ \Big\{- 2 \pi^2 T -  \big( 2/3 + o(1) \big) T^{\frac{1}{3}} \log T \Big\} },
$
as required to demonstrate that (\ref{clc}) holds.
We claim that
\begin{eqnarray}
 & & 
\mathbb{P} \Big( \bigcap_{i=0}^{ m-1}{A_i} \Big\vert \bigcap_{i=1}^{ m}{
B \Big( \frac{iT}{ m} \Big) \in q +
\phi_{i+1}  } \Big) \label{yuy} \\
 & \geq & \exp \Big\{ - \Big( \frac{2}{3} + o(1) \Big)
T^{1/3} \log T \Big\}. \nonumber
\end{eqnarray}
To see this, note that
\begin{eqnarray}
 & & \mathbb{P} \Big( \bigcap_{i=0}^{m-1}{A_i} \Big\vert \bigcap_{i=1}^{m}{
B \Big( \frac{iT}{ m} \Big) \in q +
\phi_{i+1}  } \Big) \label{prf} \\
 & = & \prod_{i=0}^{m-1} 
 \mathbb{P} \Big( {A_i} \Big\vert 
 \Big\{ B \Big( \frac{iT}{ m} \Big) \in q +
\phi_{i+1}  \Big\} \cap \Big\{ B \Big( \frac{(i+1)T}{ m} \Big) \in q +
\phi_{i+2}  \Big\}  \Big). \nonumber
\end{eqnarray}
We condition on, for example, 
\begin{equation}\label{ads}
B \Big( \frac{T}{ m} \Big) \in q + \phi_2, 
B \Big( \frac{2T}{ m} \Big) \in q + \phi_3.
\end{equation}
Consider the one-dimensional process $\hat{B} : I \to \mathbb{R}$, 
$ I := \big[ T /m, 2T / m \big]$, 
that is the component of $B$ in the direction orthogonal to the line
segment 
$\big[
B \big( T/{m} \big) , 
B \big( {2T}/{ m} \big) \big].
$
Under conditioning on (\ref{ads}), $\hat{B}$ is distributed as a
Brownian bridge run to and from two given points, each of which we may
insist lies in the interval [-2,0]. Note that the motion of $B$ will remain in the
half-plane $M_1$ during the interval of time $I$ provided that
$\hat{B} \leq 1/2$ throughout this time. The probability that the
maximum of a one-dimensional Brownian bridge that is run for a given
time, starts at $a \in \mathbb{R}$ and ends at $b \in \mathbb{R}$,
does not exceed a given value is a decreasing function of $a$ and of $b$. We
may assume therefore that $\hat{B}\big(  T / m \big)=0$ 
and  $\hat{B}\big(  2T / m \big)=0$.  

Recall that, if $M^{+}$ denotes the maximum of the
one-dimensional Brownian bridge run for time $T$, then, for any $r > 0$,
\begin{equation}\label{njh}
\mathbb{P} \big( M^{+} > r \big)  = \exp{ \Big\{ - \frac{2 r^2}{T}
\Big\} }.
\end{equation}
This assertion appears as formula (3.40) in
\cite[Chapter $4$]{karatshreve}.

We deduce that 
\begin{displaymath}
\mathbb{P} \Big( \hat{B}(t) \leq 1/2 \, \textrm{for all $t \in I$}
\Big) \geq 1 - \exp\Big\{ - \frac{m}{2T}
\Big\}
 \geq \frac{1}{8} T^{- \frac{2}{3}},
\end{displaymath}
the latter inequality being valid for high values of $T$.

Thus, given the occurrence of
(\ref{ads}), the probability that the event $A_1$ occurs is at least
$\frac{1}{8} T^{-2/3}$. From this, and the product form of
(\ref{prf}), follows (\ref{yuy}).

From (\ref{clc}) and (\ref{yuy}), we find that
$$
\mathbb{P} \Big( \bigcap_{i=0}^{m-1}{A_i} \Big) \geq 
\exp{ \Big\{ - 2 \pi^2 T - \Big( \frac{4}{3} + o(1) \Big) T^{\frac{1}{3}} \log T \Big\}}.
$$
The statement of the lemma follows from the inclusion (\ref{useit}). $\Box$ \\
The following result is required.
\begin{lemma}\label{stp}
For any constant $C_2 \geq  128 \pi$, the
planar Brownian loop $B:[0,T] \to \mathbb{R}^2$ has the property that
$$
\sup_{t \in [0,T]}{\vert B(t) \vert^2}
$$
is bounded above in distribution by $C_2 T + Z^2 $,
where $Z$ is a normal random variable
with mean zero and variance $T$. That is, for all $a \in (0,\infty)$,
$$
\mathbb{P} \Big( \sup_{t \in [0,T]}{\vert B(t) \vert^2} \geq a \Big)
\leq  \mathbb{P} \Big( C_2 T + Z^2 \geq a \Big) .
$$
\end{lemma}
{\bf Remark} In fact, the method of proof we give would permit us to
pick $Z$ to be distributed normally with mean zero and variance
$\big(1/2 + \epsilon \big) T$, for any $\epsilon>0$ (provided that the
constant $C_2$ is changed suitably.) We made the choice $\epsilon =
1/2$, because it is not valuable in the
application to choose any lower value for $\epsilon$.\\
\noindent{{\bf Proof}} Note that
\begin{equation}\label{ineq}
\sup_{t \in [0,T]}{ \big\vert B(t) \big\vert^2} \leq
\sup_{t \in [0,T]}{ B_1(t)^2} + 
\sup_{t \in [0,T]}{ B_2(t)^2},
\end{equation}
where $B_1$ and $B_2$ denote the components of $B$.
From the probability (\ref{njh}) that the maximum of a one-dimensional
Brownian bridge exceeds a given level, and its counterpart for the
minimum value reaching below a prescribed value, it follows that 
\begin{equation}\label{ct}
\mathbb{P} \Big( \sup_{t \in [0,T]}{ B_j(t)^2 } > r^2 \Big) \leq 2
\exp{ \Big\{ - \frac{2 r^2}{T} \Big\} }, \ \textrm{for $j \in \{ 1,2 \}$.}
\end{equation}
Replacing $r$ by $r/{\sqrt{2}}$ in (\ref{ct}), we obtain from (\ref{ineq}) that
\begin{equation}\label{cto}
\mathbb{P} \Big( \sup_{t \in [0,T]}{ \big\vert B(t) \big\vert^2 } > r^2 \Big) \leq 4
\exp{ \Big\{ - \frac{r^2}{T} \Big\} }.
\end{equation}
A lower bound on the tail of a random variable $Z$, distributed
normally with mean zero and variance $T$, is now obtained:
\begin{eqnarray}
\mathbb{P} \big( Z^2 > r^2 \big) & = &
\mathbb{P} \Big( \frac{\vert Z \vert}{\sqrt{T}} > \frac{r}{\sqrt{T}}
\Big) \nonumber \\
 & \geq & \frac{\frac{2r}{\sqrt{T}}}{\Big(\frac{r^2}{T} + 1\Big)\sqrt{2\pi}}\exp\Big\{-
 \frac{r^2}{2T} \Big\},  \label{poklj}
\end{eqnarray} 
where the inequality follows from a standard bound  
on the tail of the normal distribution, presented in
Section 14.8 of \cite{williams}. 
From (\ref{cto}) and (\ref{poklj}), it follows that
$$
\mathbb{P} \Big( \sup_{t \in [0,T]}{ \big\vert B(t) \big\vert^2 } > r^2 \Big)
\leq \mathbb{P} \big( Z^2 > r^2 \big),
$$
for $r \geq C \sqrt{T}$, where $C = 8 \sqrt{2\pi}$.
We have derived
\begin{equation}\label{gh}
\mathbb{P} \Big( \sup_{t \in [0,T]}{\vert B(t) \vert^2} \geq a \Big)
\leq  \mathbb{P} \Big( C_2 T + Z^2 \geq a \Big) .
\end{equation}
for any $C_2 \geq 0$ and each $a \geq C^2 T$. By choosing $C_2 = C^2$,
the right-hand-side (\ref{gh}) becomes equal to $1$ for any $a \in [0,
C^2 T]$. This establishes the statement of the Lemma. $\Box$ \\
We require some notation before proceeding.
\begin{definition}\label{defnsix}
Let $m \in \mathbb{N}$ and let $t' \in [0,T]$.
\begin{itemize}
\item Let $P = P_m^{t'}$ denote the polygon
whose vertices are given by
$$
\Big\{ B \Big( \frac{jT}{m} + t' \Big) : j \in \{0,\ldots,m-1 \} \Big\}.
$$
\item
Let the length of the edges of
$P_m^{t'}$ be denoted by $\Big\{ L_i : i \in \{ 1, \ldots, m \}
\Big\}$, so that 
$$
L_i = \Big\vert B \Big( \frac{iT}{m} + t' \Big) -  B \Big(
\frac{(i-1)T}{m} + t' \Big) \Big\vert,
$$
and let $l_i$ denote the line segment whose endpoints are $B \big(
(i-1)T/m + t' \big)$ and  $B \big( iT/m + t' \big)$.
\item
for $i \in \{ 1, \ldots, m \}$, let $R_i = R_i(t')$ denote the maximum of the
distance of the range of the motion $B$ from the segment $l_i$ during
that interval of time in which $l_i$ is traversed by $B$. This is,
$$
R_i = \sup_{t \in [0,T/m]}{d \Big( B \Big( \frac{(i-1)T}{m} + t' + t
\Big), l_i \Big)},
$$
where $d$ denotes the distance between two sets in $\mathbb{R}^2$.
\item
for $i \in \{ 1, \ldots, m \}$, let $\hat{R}_i = \hat{R}_i(t')$ denote the maximum of
the absolute value of the displacement of the motion from a point that
traverses $l_i$ at a linear rate during the time $\big[ (i-1)T/m, iT/m
\big]$. That is, 
$$
\hat{R}_i  =  \sup_{t \in [0,T/m]} d \left( B \Big( \frac{(i-1)T}{m} + t' + t
\Big),  \right.
$$ 
$$
 \left. \Big(
 1 - \frac{mt}{T} \Big)  B \Big( \frac{(i-1)T}{m} + t' \Big) +
 \frac{mt}{T} B \Big( \frac{iT}{m} + t' \Big)  
  \right) .
$$
\item let the set $Q_i = Q_i(t')$ be given by $Q_i = \{ x \in \mathbb{R}^2 : d (x,l_i) \leq R_i \}$.
\end{itemize}
\end{definition}
\begin{figure}\label{figtwo}
\begin{center}
\includegraphics{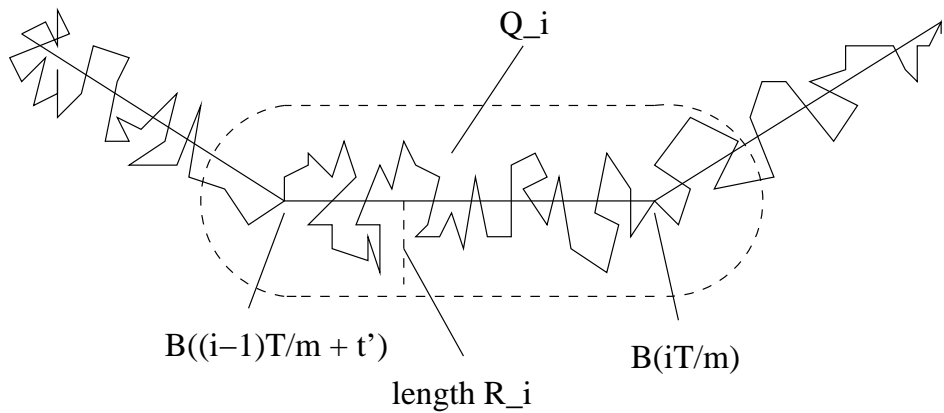}
\end{center}
\caption{The objects of Definition \ref{defnsix}}
\end{figure}
\begin{lemma}\label{assert}
For any $T \in (0,\infty), t' \in [0,T/m]$ and $m
\in \mathbb{N}$,
\begin{equation}\label{mom}
\enc(B) \subseteq \conv(P) \cup \bigcup_{i=1}^{m}{ Q_i }.
\end{equation}
\end{lemma}
{\bf Proof} The objects used in the proof are depicted in Figure
$3$. Take a point $y \in \enc(B) \setminus \conv(P)$. 
Locate a half-plane $H$ that contains $\conv(P)$ and excludes
$y$. Let $l'$ denote the line through $y$ that intersects the
boundary of $H$ at right angles. Note that there exists a point $z$ in
the range of $B$ that lies on $l'$, at a distance from the half-plane
$H$ greater than that of $y$: this is because $y \in \enc(B)$. Let
$l_i$ be the interval in the polygon $P$ traversed while the Brownian
path captures the point $z$: that is $i \in \{1,\ldots,m \}$ is chosen
so that 
\begin{equation}\label{wen}
z = B(t), \ \textrm{where $t \in \Big[ \frac{(i-1)T}{m} + t' ,
\frac{iT}{m} + t' \Big]$.}
\end{equation}
\begin{figure}\label{figthr}
\begin{center}
\includegraphics{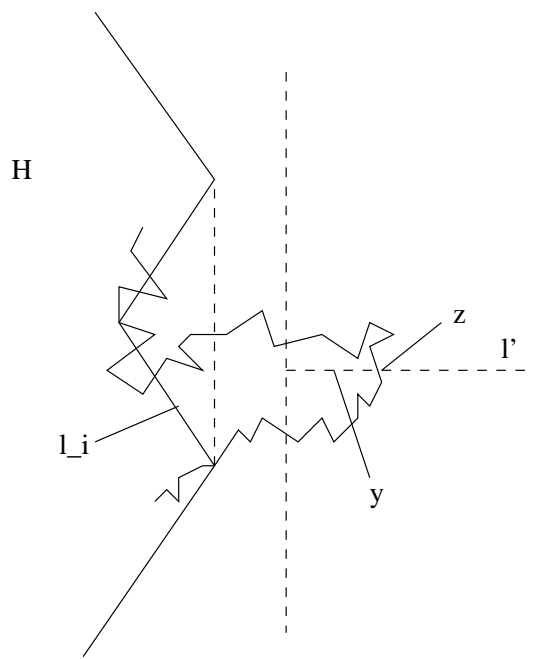}
\end{center}
\caption{The proof of Lemma \ref{assert}}
\end{figure} 
The point $y$ is closer to any given point in $H$ than is $z$, and
each point in the line segment $l_i$ lies in $\conv (P)$, 
and so in $H$. Hence, $d(y,l_i) \leq d(z,l_i)$. That $d(z,l_i) \leq R_i$
follows from (\ref{wen}). Hence, $y \in Q_i$. $\Box$ \\  
\begin{lemma}\label{lemcl}
Let $\epsilon \in [0,1/6)$, and let $m = \lfloor T^{1/3 - \epsilon} \rfloor$. Suppose that the point
$t' \in [0,T]$ in the
definition of the polygon $P$ is fixed, or that it is sampled
randomly, according to an arbitrary ditribution on $[0,T]$,
and independently of the randomness that generates $B:[0,T] \to \mathbb{R}^2$.
For $\epsilon' \in (0,1/3)$, let $H = H(\epsilon')$ denote the event that 
$$
\big\vert \conv(P) \big\vert \geq \pi T^2 - T^{4/3 + \epsilon'},
$$ 
where $\conv(P)$ denotes the convex hull of the polygon
$P$. Provided that $\epsilon' > 2 \epsilon$, for sufficiently
large values of $T$,
$$
\mathbb{P} \Big( H(\epsilon')^{c} \cap \big\{ \vert \enc(B) \vert \geq
\pi T^2 \big\} \Big) \leq \exp \Big\{- 2 \pi^2 T - c_1 T^{\frac{1}{3} + 2(\epsilon' -
\epsilon)} \Big\},
$$
where $c_1$ is any positive constant at most $\frac{2^{-10}}{108 \cdot
(128)^2 \pi^4}$. 
\end{lemma}
{\bf Proof}
For $a \in \mathbb{N}$, let $S_a$ denote the event that
$$
\pi T^2 - 2^a T^{4/3 + \epsilon'} > \big\vert \conv(P) \big\vert \geq
\pi T^2 - 2^{a+1} T^{4/3 + \epsilon'}.
$$
Note that 
\begin{equation}\label{ahtch}
H^c = \Big( \bigcup_{a = 0}^{k}{S_a} \Big) \cup R,
\end{equation}
where $k = \lfloor \log_2 {T^{2/3 - \epsilon'}} \rfloor
- 1$, with $R$ denoting the event that the area of $\conv (P)$ 
is at most $\big( \pi - 1/2 \big)T^2$.

Note that 
\begin{eqnarray}
 \Big\vert \, \enc(B) \setminus \conv(P) \Big\vert & \leq & \sum_{i=1}^{m}{\big\vert
 Q_i \big\vert}
 \nonumber \\
 & \leq & 2 \sum_{i=1}^{m}{ \big( L_i + 2 R_i \big) R_i} \label{klo} \\
 & \leq & 2 \sum_{i=1}^{m}{ \big( L_i + 2 \hat{R}_i \big) \hat{R}_i} \nonumber \\
 & \leq & 2 \Big(  \sum_{i=1}^{m}{L_i^2} \Big)^{1/2}  \Big(
 \sum_{i=1}^{m}{\hat{R}_i^2} \Big)^{1/2} + 4 \sum_{i=1}^{m}{\hat{R}_i^2}, \nonumber
\end{eqnarray}
where the first inequality follows from Lemma  \ref{assert}. 
The second inequality follows from the fact that $Q_i$ is contained in a rectangle of length
$L_i + 2 R_i$ and width $2 R_i$, while the third is implied by the
inequality $R_i \leq \hat{R}_i$.

Let $a \in \big\{0,\ldots,k \big\}$.
Note that, on the event $S_a
\cap \{ \vert \enc(B) \vert \geq \pi T^2 \}$, we have that
\begin{equation}\label{ewq}
\Big\vert \enc(B) \setminus  \conv (P) \Big\vert \geq \big\vert \enc(B) \big\vert -
\big\vert \conv(P) \big\vert \geq 2^a T^{4/3 + \epsilon'}.
\end{equation}
We may write
\begin{equation}\label{rew}
L_i^2  =  \frac{T \big( E_{2i -1}^2 + E_{2i}^2 \big)}{m},
\end{equation}
where, for each $i \in \{ 1, \ldots, m \}$, 
the quantities $\sqrt{T/m} E_{2i-1}$ and $\sqrt{T/m} E_{2i}$ are the horizontal and vertical
components of the vector $B\big(iT/m + t'\big) - B\big((i-1)T/m +
t'\big)$. As such, the family
$\big\{ E_i : i \in \{ 1, \ldots, 2m \} \big\}$ is a
collection of independent standard normal random variables,
conditioned by insisting that 
\begin{equation}\label{ehalf}
\sum_{i=1}^{m}{E_{2i - 1}} = 0 , \, \,
\sum_{i=1}^{m}{E_{2i}} = 0. 
\end{equation}
For each $i \in \{1,\ldots,m\}$ and any pair $(x,y) \in \mathbb{R}^2$,
the conditional distribution of the random variable $\hat{R}_i$ given
the event that $E_{2i-1} = x$ and $E_{2i} = y$ is independent of $x$
and $y$. Indeed, the process $Z_i: [(i-1)T/m,iT/m] \to \mathbb{R}^2$
given by 
\[
 Z_i(t) = B(t) + 
          \Big(1 -  \big( mt/T - (i-1) \big) \Big) B \Big(
	  \frac{(i-1)T}{m} \Big) +  \big(mt/T - (i-1)  \big)
            B \Big( \frac{iT}{m} \Big) 
\]
is a standard Brownian loop run for a time of $T/m$, no matter how we condition the values of
the endpoints $B\big((i-1)T/m \big)$ and $B\big( iT/m \big)$. 
This means that, under any conditioning of the form  $E_{2i-1} = x$ and $E_{2i} = y$, each $\hat{R}_i$ has the distribution of
the maximal
Euclidean distance of this Brownian loop.
As such, we may apply Lemma \ref{stp}
to find a collection of standard normal random variables
$\big\{ F_i : i \in \{ 1, \ldots, m \} \big\}$ for which  
\begin{equation}\label{req} 
 \hat{R}_i^2  \leq   \frac{T}{m} \big( C_2 + F_i^2  \big).
\end{equation}
The fact that the conditional distribution of
$\hat{R}_i$ does not depend on $x$ and $y$ implies that we may assume that each $F_i$ is independent
of $F_j$, for $j \in \{1,\ldots,m\}$ with $j \not= i$, and of each
$E_j$, for $j \in \{1,\ldots,2m \}$.

From (\ref{klo}), (\ref{ewq}), (\ref{rew}) and (\ref{req}) follows
$$
\frac{2T}{m} \Big( \sum_{i=1}^{2m}{E_i^2} \Big)^{1/2}
 \Big( C_2 m + \sum_{i=1}^{m}{F_i^2} \Big)^{1/2} + 4 C_2 T + \frac{4T}{m}
 \sum_{i=1}^{m}{F_i^2}  \geq 2^a T^{4/3 + \epsilon'}. 
$$
Using the inequalities $m \geq T^{1/3 - \epsilon}/2$ and $C_2 \geq
\sqrt{C_2} \geq 1$,
\begin{equation}\label{asd}
  \Big( \sum_{i=1}^{2m}{E_i^2} \Big)^{1/2}
 \Big( m + \sum_{i=1}^{m}{F_i^2} \Big)^{1/2} + 2 m +
 2 \sum_{i=1}^{m}{F_i^2} \geq \frac{2^{a-2} T^{2/3 - \epsilon +
 \epsilon'}}{C_2}.
\end{equation}
We note that 
\begin{equation}\label{cbv}
\sum_{i=1}^{m}{L_i^2} \geq \frac{\big( \sum_{i=1}^{m}{L_i} \big)^2}{m}
\geq  \frac{\Big( \arcl \big( \conv(P) \big) \Big)^2}{m} \geq \frac{ 4
\pi \big\vert \conv(P) \big\vert}{m} ,
\end{equation}
where $\arcl$ denotes arclength.
The successive inequalities in (\ref{cbv}) are consequences of the
Cauchy-Schwarz inequality, the fact that the arclength of 
any polygon
exceeds that of its convex hull, and the standard isoperimetric
inequality, which in this case asserts that
$$
 \Big( \arcl \big( \conv P \big) \Big)^2 \geq 4 \pi \big\vert 
 \conv(P) \big\vert.
$$  
It follows from (\ref{rew}),(\ref{cbv}) and the lower bound on $\vert \conv(P)
\vert$ provided by the occurrence of the event $S_a$ that
\begin{equation}\label{asf}
\sum_{i=1}^{2m}{E_i^2} \geq  
4 \pi^2 T - 2^{a+3}  \pi T^{1/3 + \epsilon'},
\end{equation}
because the left and right-hand-side of this inequality are
respectively bounded below and above by the quantity $\frac{ 4 \pi
\big\vert  \conv(P) \big\vert}{T}$.
We see that the event 
$$
\Big\{ \vert enc(B) \vert \geq \pi T^2 \Big\} \cap S_a
$$
is contained in the event $F$, specified by the occurrence of the
inequalities in (\ref{asd}) and (\ref{asf}). Note that, we have the following
inclusion:
\begin{equation}\label{finc}
F \subseteq A_C^1 \cup A_C^2 \cup A_C^3,
\end{equation}
where the events on the right-hand-side are given by
\begin{eqnarray}
A_C^1 & = & \Big\{ \sum_{i=1}^{2m}{E_i^2} \geq C T \Big\}, \nonumber \\
A_C^2 & = &  \Big\{ \sum_{i=1}^{2m}{E_i^2} \in \big[ 4 \pi^2 T - 2^{a+3} \pi
T^{1/3 + \epsilon'} , C T \big) \Big\} \nonumber \\
 & & \cap \ \Big\{ m + \sum_{i=1}^{m}{F_i^2} \geq  \sum_{i=1}^{2m}{E_i^2}
 \Big\} \nonumber 
\end{eqnarray}
and 
\begin{eqnarray}
A_C^3 & = &   \Big\{ \sum_{i=1}^{2m}{E_i^2} \in \big[ 4 \pi^2 T - 2^{a+3} \pi
T^{1/3 + \epsilon'} , C T \big) \Big\} \nonumber \\
 & & \cap \ \Big\{ m + \sum_{i=1}^{m}{F_i^2} \geq  \frac{2^{2a-4} T^{1/3 - 2
 \epsilon + 2 \epsilon'}}{9 C C_2^2 } \Big\}. \nonumber 
\end{eqnarray}
Here, $C$ denotes a fixed constant, and the constant $C_2$ satisfies
the bound stated in Lemma \ref{stp}. To derive (\ref{finc}), note that,
if the event $F \cap \big( A_C^1 \cup A_C^2 \big)^c$ occurs, then
\begin{eqnarray}
 & & 3 \Big( \sum_{i=1}^{2m}{E_i^2} \Big)^{\frac{1}{2}}  \Big( m +
   \sum_{i=1}^{m}{F_i^2}  \Big)^{\frac{1}{2}} \nonumber \\
 & \geq &  \Big( \sum_{i=1}^{2m}{E_i^2} \Big)^{\frac{1}{2}}  \Big( m +
   \sum_{i=1}^{m}{F_i^2}  \Big)^{\frac{1}{2}} + 2  
     \Big( m + \sum_{i=1}^{m}{F_i^2}  \Big)  \nonumber \\
   & \geq &   \frac{2^{a-2}}{C_2} T^{\frac{2}{3} - \epsilon +
   \epsilon'}. \label{isidkaal} 
\end{eqnarray}
The first inequality follows from the fact that 
 $m + \sum_{i=1}^{m}{F_i^2} < \sum_{i=1}^{2m}{E_i^2}$, which holds on
 the event  $F \cap \big( A_C^1 \cup A_C^2 \big)^c$, while the second is
 simply (\ref{asd}). Using (\ref{isidkaal}), and the fact that $ \sum_{i=1}^{2m}{E_i^2} <
 CT$, we find that the second event of the two whose intersection
 defines $A_C^3$ must occur. The first occurs provided that both the event
 $\big( A_C^1 \big)^c$ occurs and the inequality (\ref{asf}) is satisfied, which means
 that it occurs if the event  $F \cap \big( A_C^1 \cup A_C^2 \big)^c$ does.

The standard normal random variables $\big\{ E_i: i \in
\{ 1, \ldots, 2m \} \big\}$ are independent, conditioned only by the
two linear constraints (\ref{ehalf}). Thus, the quantity $
\sum_{i=1}^{2m}{E_i^2}$ has the $\chi^2$-distribution with $2m - 2$
degrees of freedom. The random variable $\sum_{i=1}^{m}{F_i^2}$ has
the $\chi^2$-distribution with $m$
degrees of freedom. Recall that the density $f_m : [0,\infty) \to
[0,\infty)$ of this latter distribution is given by
\begin{equation}\label{fchi}
f_m \big( x \big) = \frac{x^{m/2 - 1} \exp \big\{ -x/2 \big\} }{2^{m/2}
\Gamma \big( m/2 \big) }. 
\end{equation}
From this formula, we find that, for any given $C>0$ and all
sufficiently high values of $T$,
\begin{equation}\label{nal}
\mathbb{P}\big( A_C^1 \big) \leq \exp \Big\{ - \frac{CT}{4} \Big\}.
\end{equation}
Recalling the independence property of the family of random variables
$\big\{ E_i : i \in \{1,\ldots,2m\} \big\}$ and $\big\{ F_i : i \in
\{1,\ldots,m\} \big\}$ asserted after (\ref{req}), it follows that
\begin{equation}\label{aytwo}
\mathbb{P} \big( A_C^2 \big) \leq \int_{l}^{CT}{f_{2m -2}(x) dx}  \int_{l - m}^{\infty}{f_m(x) dx},
\end{equation}
where $l = 4 \pi^2 T - 2^{a+3} \pi T^{1/3 + \epsilon'}$. 
The right-hand-side of (\ref{aytwo}) may be bounded above by
\begin{eqnarray}
 & & \! \! \! \! \! \! \! \! \! \! CT \big( CT \big)^{T^{1/3 - \epsilon}}  \exp{\Big\{ - 2 \pi^2 T +
 2^{a+2} \pi T^{\frac{1}{3} + \epsilon' } \Big\} } \, \, \times \nonumber \\
 & & \! \! \! \!  \! \! \! \!  \! \! \Big[  CT \big( CT \big)^{T^{1/3 - \epsilon}}  \exp{\Big\{ - 2 \pi^2 T +
 2^{a+2} \pi T^{\frac{1}{3} + \epsilon' } + \frac{1}{2} T^{\frac{1}{3} - \epsilon}
 \Big\} }  + \exp \Big\{-CT/4 \Big\} \, \Big] \nonumber
\end{eqnarray}
by making use of $m \leq T^{1/3 - \epsilon}$.
From the fact that $2^{a+1} \leq T^{2/3 - \epsilon'}$ follows
that for any constant $c'$ and for all $T$ exceeding a $c'$-dependent
constant, 
\begin{equation}\label{nbl}
\mathbb{P} \big( A_C^2 \big) \leq   \exp \Big\{ - \Big[ 4 \pi \big( \pi - 1\big) - c' \Big] T \Big\},
\end{equation}
provided that $C > 8 \pi \big( \pi - 1 \big)$.
On the event $A_3$, $\sum_{i=1}^{m}{F_i^2} \geq m$: this follows from
the second requirement in the definition of $A_3$, along with the
weaker than given inequality $\epsilon < 2 \epsilon'$, and $m \leq
T^{1/3 - \epsilon}$. Thus, if $A_3$ occurs, then 
$$
 \sum_{i=1}^{m}{F_i^2} \geq
 \frac{2^{2a-5} T^{1/3 - 2
 \epsilon + 2 \epsilon'}}{9 C C_2^2 }.
$$
Thus, by a similar estimate as in the case of $A_C^2$, 
\begin{eqnarray}
\mathbb{P} \big( A_C^3 \big)
 & \leq & CT \big( CT \big)^{T^{1/3 - \epsilon}}  \exp{\Big\{ - 2 \pi^2 T +
 2^{a+2} \pi T^{1/3 + \epsilon' } \Big\} } \nonumber \\
 & & \times \Big[  CT \big( CT \big)^{T^{1/3 - \epsilon}}  \exp{\Big\{
 - \frac{2^{2a - 6}}{9 C C_2^2} T^{\frac{1}{3} - 2 \epsilon + 2 \epsilon'}
 \Big\} }  + \exp\Big\{ -CT/4 \Big\} \Big]. \nonumber
\end{eqnarray}
Using the fact that $\epsilon' > 2\epsilon$,
\begin{equation}\label{naal}
\mathbb{P} \big( A_C^3 \big) \leq \exp \Big\{ - 2 \pi^2 T  - \frac{2^{2a - 7}}{9 C C_2^2} T^{\frac{1}{3} - 2
 \epsilon + 2 \epsilon'} \Big\} \, + \, \exp \Big\{  - CT/5
 \Big\} .
\end{equation}
We find that 
\begin{eqnarray}
 & & \mathbb{P} \Big( S_a \cap \Big\{ \big\vert \enc(B) \big\vert \geq
 \pi T^2 \Big\} \Big) \nonumber \\
 & \leq & \mathbb{P} \Big( A_C^1 \cup A_C^2 \cup A_C^3
 \Big) \leq  \mathbb{P} \big( A_C^1 \big) +  \mathbb{P} \big( A_C^2 \big)
 + \mathbb{P} \big( A_C^3 \big) \nonumber \\
 & \leq & \exp \Big\{ - 2 \pi^2 T  - \frac{2^{2a - 8}}{9 C C_2^2} T^{\frac{1}{3} - 2
 \epsilon + 2 \epsilon'} \Big\}  \, + \, \exp \Big\{ - CT/6 \Big\} \nonumber 
\end{eqnarray}
where in the first inequality, we used (\ref{finc}), and in the third,
(\ref{nal}), (\ref{nbl}) and (\ref{naal}).  
Recalling that $k = \lfloor \log_2 T^{2/3 - \epsilon'}
 \rfloor - 1$, we find that, provided that $C > 12 \pi^2$, 
\begin{displaymath}
 \mathbb{P} \bigg( \Big( \bigcup_{a=0}^{k} S_a \Big) \cap \Big\{ \big\vert \enc(B) \big\vert \geq
 \pi T^2 \Big\} \bigg) \leq  \exp \Big\{ - 2 \pi^2 T  - c_1 T^{\frac{1}{3} - 2
 \epsilon + 2 \epsilon'} \Big\}, 
\end{displaymath}
for any constant $c_1$ satisfying 
$$
c_1 \in \bigg( 0, \frac{2^{-9}}{9 C C_2^2}
\bigg) \subseteq \bigg(0,  \frac{2^{-10}}{9 \cdot 12 \cdot
(128)^2 \pi^4  } \bigg), $$
where the second inclusion is ensured by choosing $C_2 = 128 \pi$ and $C$ slightly above $12 \pi^2$. 

By (\ref{ahtch}), it only remains
to bound the probability of the event $R$ that appears in
there. Note that if $R$ occurs, then
$$
\big\vert \enc(B) \setminus \conv(P) \big\vert \geq \frac{1}{2} T^2.
$$
From (\ref{klo}), 
we see that the occurrence of $R$ implies that
$$
 2 \Big(  \sum_{i=1}^{m}{L_i^2} \Big)^{1/2}  \Big(
 \sum_{i=1}^{m}{\hat{R}_i^2} \Big)^{1/2} + 4
 \sum_{i=1}^{m}{\hat{R}_i^2} \geq \frac{1}{2} T^2,
$$
and so that
\begin{equation}\label{dhtrifl}
  \Big(  \sum_{i=1}^{2m}{E_i^2} \Big)^{1/2}  \Big( C_2 m + 
 \sum_{i=1}^{m}{F_i^2} \Big)^{1/2} + 2 \Big( C_2 m +
 \sum_{i=1}^{m}{F_i^2} \Big) \geq
  \frac{1}{8} T^{\frac{4}{3} - \epsilon}, 
\end{equation}
by use of (\ref{rew}), (\ref{req}) and the inequality $m \geq \frac{1}{2} T^{1/3 - \epsilon}$.
If the inequality (\ref{dhtrifl}) is satisfied, then so is one of 
\begin{displaymath}
C_2 m + \sum_{i=1}^{m}{F_i^2}  \geq  \frac{1}{24} T^{\frac{4}{3} -
\epsilon} 
\end{displaymath}
and
\begin{displaymath}
\sum_{i=1}^{2m}{E_i^2}  \geq  \frac{1}{24} T^{\frac{4}{3} -
\epsilon}.
\end{displaymath}
Since $m \leq T^{1/3}$ and $\epsilon < 1/3$, each of these inequalitities is
satisfied with a probability that decays at a rate faster than exponential.
This completes the proof.
 $\Box$ \\
It is of interest to bound the extent of the excess of area
captured by the motion, partly for the reasons presented in the
heuristic discussion that ends the Introduction. In addition, we will
make use of the
following bound in the proof of Theorem \ref{tgg}. 
\begin{prop}\label{propexcarea}
For $\delta \in (0,1/3)$
and for all $T \in (0,\infty)$ 
sufficiently high, we have that
\begin{displaymath}
 \mathbb{P} \Big( \big\vert \enc (B) \big\vert \geq \pi T^2 +
 T^{\frac{4}{3} + \delta} \Big) \leq \exp \Big\{ - 2 \pi^2 T - \frac{\pi}{4}
 T^{\frac{1}{3} + \delta} \Big\}.
\end{displaymath}
\end{prop}
{\bf Proof}
We set $m = \lfloor T^{\frac{1}{3}} \rfloor$.
Note that the event whose probability we seek to bound lies in $Y_1
\cup Y_2$, where the events $Y_1$ and
$Y_2$ are given by
\[
 Y_1 =  \Big\{ \big\vert \enc (B) \big\vert \geq \pi T^2 +
 T^{\frac{4}{3} + \delta} \Big\} \cap  \Big\{ \big\vert \conv (P)
 \big\vert < \pi T^2 +
 \frac{1}{2} T^{\frac{4}{3} + \delta} \Big\},
\]
and by
\[
Y_2 =  \Big\{ \big\vert \conv (P)
 \big\vert \geq \pi T^2 +
 \frac{1}{2} T^{\frac{4}{3} + \delta} \Big\}.
\]
Let $\hat{T} \in (0,\infty)$ be given by $\pi \hat{T}^2 = \pi T^2 +
T^{\frac{4}{3} + \delta}$. We find that
\begin{eqnarray}
 \mathbb{P} (Y_1)  & \leq & 
 \mathbb{P} \Big(  \Big\{ \big\vert \enc (B) \big\vert \geq \pi \hat{T}^2 \Big\} \cap  \Big\{ \big\vert \conv (P)
 \big\vert < \pi \hat{T}^2 -
 \frac{1}{4} \hat{T}^{\frac{4}{3} + \delta} \Big\} \Big)  \label{jjasp} \\ 
& \leq & \exp \Big\{ - 2 \pi^2 T - \frac{c_1}{4^2}
 T^{\frac{1}{3} + 2 \delta} \Big\}. \nonumber
\end{eqnarray}
The first inequality in (\ref{jjasp}) 
follows for high values of $T$ from the fact that
$\delta < 2/3$. The second is an application of Lemma \ref{lemcl} (with
the choice $\epsilon'  = \delta - \frac{\log 4}{\log \hat{T}}$), and of
the inequality $\hat{T} \geq T$. 

Note that, by (\ref{rew}) and (\ref{cbv}),
\[
 \sum_{i=1}^{2m}{E_i^2} \geq \frac{4\pi}{T} \vert \conv(P) \vert,
\]
so that
\begin{equation}\label{etrt}
 \mathbb{P} \big( Y_2 \big) 
 \leq   \mathbb{P} \Big( \sum_{i=1}^{2m}{E_i^2} \geq 4 \pi^2 T +
 2 \pi T^{\frac{1}{3} + \delta} \Big)
  \leq  \exp \Big\{ - 2 \pi^2 T -  \frac{\pi}{2}
 T^{\frac{1}{3} + \delta} \Big\}. 
\end{equation}
We used the form (\ref{fchi}) of the density of the $\chi^2$-distribution
in the second inequality, as well as the fact that $m \leq T^{1/3}$.
Applying the bounds (\ref{jjasp}) and (\ref{etrt}) to estimate
$\mathbb{P}(Y_1\cup Y_2)$,
we deduce the statement of the proposition. $\Box$ \\
\begin{lemma}\label{remaftdef}
Given functions $f,g : [0,T] \to [0,\infty)$, let $Q_{f,g}$ denote the
event that 
\begin{displaymath}
\sup_{t \in [0,T]}{\sup_{ s \in [t,t+f(T)]} } \big\vert B(s)
- B(t) \big\vert \leq g (T).
\end{displaymath}
Provided that $f(T) \leq T$, we have that 
\begin{equation}
\mathbb{P} \big( Q_{f,g}^c \big) \leq \bigg( \frac{32 \sqrt{2} T}{\sqrt{\pi
f(T)} g(T)} + \frac{16\sqrt{T}}{\sqrt{\pi}g(T)} \bigg)\exp \bigg\{ - 
\frac{g(T)^2}{128 f(T)} \bigg\}.
\end{equation}
\end{lemma}
\noindent{{\bf Proof:}}
Firstly, note that
\begin{equation}\label{netove}
 Q_{f,g}^c \subseteq \bigcup_{j=1}^{\lfloor T/f(T) \rfloor} C_j ,
\end{equation}
where $C_j$ is the event that there exist $t_1,t_2 \in \big[ j f(T),
(j+2) f(T) \big]$ for which $\big\vert B(t_1) - B(t_2) \big\vert >
g(T)$.
Noting that the process $B_t:[0,T] \to \mathbb{R}^2$ given by $B_t(s)
= B(t+s) - B(t)$ has the same law as $B = B_0$, we see that 
\begin{equation}\label{epc}
\textrm{each of
the events $C_j$ for $j \in \big\{1,\ldots,\lfloor T/f(T) \rfloor
\big\}$
has an equal probability.}
\end{equation}
We represent $B:[0,T] \to \mathbb{R}^2$
in the form $B(t) = W(t) - tW(T)/T$, where $W: [0,T] \to \mathbb{R}^2$
is distributed as a standard planar Brownian motion.   Note then the inclusion 
\begin{equation}\label{netover}
C_1 \subseteq \big(  A_1 \cap F^c \big) \cup F ,
\end{equation}
where  $A_1$ is the event that there exist $t_1,t_2 \in \big[ f(T),
 3 f(T) \big]$ for which $\big\vert W(t_1) - W(t_2) \big\vert >
 g(T)/2$, and where $F$ is the event
$$
 F = \Big\{ W(T) > \frac{T g(T)}{4 f(T)} \Big\}.
$$
Indeed, on the event $C_1 \cap F^c$, we find that
\begin{eqnarray}
 \big\vert W(t_1) - W(t_2) \big\vert & \geq & \big\vert B(t_1) - B(t_2) \big\vert
 - \frac{\big\vert t_1 - t_2 \big\vert}{T} \big\vert W(T) \big\vert
 \nonumber \\ 
 & \geq & g(T) -  \big\vert t_1 - t_2 \big\vert \frac{g(T)}{4 f(T)} \geq
 \frac{1}{2} g(T), \nonumber  
\end{eqnarray} 
where $t_1,t_2 \in [0,T]$ are as in the definition of the event $C_j$.
Thus, $A_j$ occurs for that $j$ for which $t_1,t_2 \in  \big[ j f(T),
(j+2) f(T) \big]$. We have verified (\ref{netover}).

Given that, for each $j \in \big\{ 1, \ldots, \lfloor T/f(T) \rfloor \big\}$,
\begin{equation}\label{kata}
 \mathbb{P} \big( C_j \big) = \mathbb{P} \big( C_1 \big)
 \leq \mathbb{P} \big( A_1 \big) + \mathbb{P}\big( F\big),
\end{equation}
we seek to bound the probability of the events $A_1$ and $F$. To this end, 
note that, if $A_1$ occurs,
then one of the components in the $x$ or $y$ direction of
$W(t_1)-W(t_2)$ on $\big[ f(T), 3 f(T) \big]$ has a difference
between its maximum and minimum values that exceeds $g(T)/(2
\sqrt{2})$. For this component, one of the maximum and the absolute value of the minimum
exceeds $g(T)/(4\sqrt{2})$. Thus,
\begin{eqnarray}
\mathbb{P} \big( A_1 \big) & \leq & 4 \mathbb{P} \Big( \sup_{t \in
[0,2f(T)]}W_1(t) > \frac{g(T)}{4\sqrt{2}} \Big) \nonumber \\
 & = & 8 \mathbb{P} \Big( W_1 \big( 2f(T) \big) >
 \frac{g(T)}{4\sqrt{2}} \Big) \nonumber \\
 & \leq  & \frac{32 \sqrt{2} \sqrt{f(T)}}{\sqrt{\pi} g(T)} 
 \exp \Big\{ - \frac{g(T)^2}{128 f(T)}\Big\}, \label{polp} 
\end{eqnarray}
where in the equality, we used the reflection principle, and, in the
latter inequality, a standard tail bound for a normal random variable
(Theorem $1.4$ of \cite{durr}).

To bound the probability of the event $F$, 
note that if  $ \frac{f(T)}{T} \big\vert W(T) \big\vert > \frac{1}{4}
g(T) $, then at least one of the inequalities
\begin{displaymath}
  \frac{f(T)}{T} \big\vert W_i (T) \big\vert > \frac{1}{4\sqrt{2}}
g(T).
\end{displaymath}
holds for $i \in \big\{ 1,2 \big\}$. By the same tail bound,
\begin{eqnarray}
\mathbb{P} \big( F \big) & = & \mathbb{P} \Big(  \frac{f(T)}{T} \big\vert W(T) \big\vert > \frac{1}{4}
g(T) \Big) \label{secun} \\
 & \leq  & \frac{16 f(T)}{\sqrt{\pi T} g(T)} \exp \Big\{ - \frac{T
g(T)^2}{64 f(T)^2} \Big\} \leq   \frac{16 f(T)}{\sqrt{\pi T} g(T)} \exp \Big\{ -
 \frac{g(T)^2}{64 f(T)} \Big\}, \nonumber
\end{eqnarray}
since $f(T) \leq T$.

We now find, as required, that
\begin{eqnarray}
 \mathbb{P} \big( Q_{f,g}^c \big) & \leq & \frac{T}{f(T)} \mathbb{P}
 \big(C_1 \big) \leq  \frac{T}{f(T)} \Big( \mathbb{P}(A_1) +
 \mathbb{P}(F) \Big) \nonumber \\
 & \leq & \bigg( \frac{32\sqrt{2}T}{\sqrt{\pi}\sqrt{f(T)}g(T)} +
 \frac{16\sqrt{T}}{\sqrt{\pi} g(T)} \bigg) \exp \Big\{ -
 \frac{g(T)^2}{128f(T)} \Big\}, \nonumber
\end{eqnarray}
by means of (\ref{netove}), (\ref{epc}) in the first inequality,
(\ref{kata}) in the second, and (\ref{polp}) and (\ref{secun}) in the
third. $\Box$
\begin{prop}\label{ttwo}
For $\epsilon_0 \in (0,1/3)$, let $\overline{J}_{\epsilon_0}$ denote the event
that, for some $\epsilon \in (\epsilon_0,1/3)$ and for some $t \in [0,T]$,
$$
\Big\vert B(t + T^{2/3 + \epsilon}) - B(t) \Big\vert \geq T^{\frac{2}{3} + 2 \epsilon}.
$$
Then
$$
\mathbb{P} \Big( \overline{J}_{\epsilon_0} \cap \Big\{ \big\vert \enc(B) \big\vert \geq
\pi T^2 \Big\}  \Big) \leq C_3 T^2 \log T \exp \Big\{  - 2 \pi^2
T - 2^{-9} T^{\frac{2}{3} + 3 \epsilon_0} \Big\},
$$ 
where $C_3$ is any constant exceeding $256 \pi^2$.
\end{prop}
\begin{lemma}\label{smalllemma}
For each $\epsilon_0 \in (0,1/3)$ and $\delta \in
 ( 0,\epsilon_0 )$, the occurrence of $\overline{J}_{\epsilon_0
 - \delta}^c$ implies that,
for each $\epsilon > \epsilon_0$, if $t,t' \in [0,T]$ satisfy
$\big\vert B(t') - B(t) \big\vert \geq \frac{9}{2} T^{2/3 + 2 \epsilon}$, then
$\big\vert t' - t \big\vert
\geq \frac{1}{2 } T^{2/3 + \epsilon}$. 
\end{lemma}
\noindent{{\bf Proof}}
Given $\epsilon, t$
and $t'$ satisfying the hypotheses of the lemma, we find that either
\begin{equation}\label{caone}
 \Big\vert B \Big( t + T^{\frac{2}{3} + \epsilon} \Big) - B
 \big( t \big) \Big\vert \geq \frac{9}{4} T^{\frac{2}{3} + 2 \epsilon},
\end{equation}
or
\begin{equation}\label{catwo}
 \Big\vert B \Big( t + T^{\frac{2}{3} + \epsilon } \Big) - B
 \big( t' \big) \Big\vert  \geq \frac{9}{4} T^{\frac{2}{3} + 2 \epsilon}.
\end{equation}
If (\ref{caone}) applies, then certainly $\overline{J}_{\epsilon_0}$
occurs, since $\epsilon > \epsilon_0$. This means of course that
$\overline{J}_{\epsilon_0  - \delta}$ occurs as well. 
Supposing the second eventuality (\ref{catwo}) and that  $\big\vert t' - t
 \big\vert < \frac{1}{2} T^{2/3 + \epsilon }$, we note that
 (\ref{catwo}) may be rewritten
\begin{equation}\label{cnb}
\Big\vert B \Big( t' + T^{\frac{2}{3} + \epsilon' } \Big) - B
\big( t' \big) \Big\vert \geq \frac{9}{4} T^{\frac{2}{3} + 2 \epsilon},
\end{equation}
where the quantity $\epsilon'$ is easily shown to satisfy the bounds
\begin{equation}\label{epsbd}
\epsilon' \in \big( \epsilon - \log 2 / \log{T}, \epsilon + 
\log(3/2) / \log{T} \big) 
\end{equation}
and thus $\epsilon' >
\epsilon - \delta > \epsilon_0 - \delta$, for high values of $T$.
The right-hand-side of (\ref{cnb}) exceeds $T^{2/3 + 2\epsilon'}$ for
high values of $T$, because (\ref{epsbd}) implies that $T^{2(\epsilon -
\epsilon')} \geq 4/9$. We have shown that the occurrence of
$\overline{J}_{\epsilon_0 - \delta}$ is a consequence of (\ref{catwo}) and the
inequality $\big\vert t' - t \big\vert < \frac{1}{2} T^{2/3 +
\epsilon }$. This establishes the statement of the lemma. $\Box$ \\
The proof of Proposition \ref{ttwo} depends largely on the following Lemma.
\begin{lemma}\label{tone}
Let $\epsilon \in ( 0, 1/3 )$. Let $J = J(\epsilon)$ denote the event that, for
some $t \in [0,T]$,
$$
\Big\vert B(t + T^{2/3 + \epsilon}) - B(t) \Big\vert \geq \frac{1}{2} T^{\frac{2}{3} + 2 \epsilon}.
$$
Then
$$
\mathbb{P} \Big( J \cap \Big\{  \big\vert \enc(B) \big\vert \geq \pi T^2
\Big\} \Big) \leq \exp \Big\{ -  2 \pi^2 T - 2^{-9}
T^{\frac{2}{3}  + 3 \epsilon} \Big\}.
$$ 
\end{lemma}
{\bf Proof}
Suppose that $\epsilon < 1/5$ (the other case is simple, and will be
handled at the end of the proof). Find $\epsilon' > 0$ satisfying 
$\epsilon' < 2/3$, $\epsilon' - 3 \epsilon < 1/3$ and $5 \epsilon +
1/3 < 2 \epsilon'$ (these conditions may be satisfied, since $\epsilon
< 1/5$). We will consider
the polygonal approximation of the range of $B$ in the case where $m =
\lfloor T^{1/3 - \epsilon} \rfloor$. In this context, the quantity
$t'$ of Definition \ref{defnsix} will be chosen to be a uniform random variable on the interval
$\big[ 0, T/m \big]$ that is independent of the motion $B$. Note that
\begin{equation}\label{hcon}
 H \big( \epsilon' \big) \subseteq \Big\{ \sum_{i=1}^{m}{L_i} \geq 2
 \pi T - \frac{\pi}{3} T^{\frac{1}{3} + \epsilon'} \Big\}.
\end{equation}
Indeed,
\begin{eqnarray}
 \sum_{i=1}^{m}{L_i} & \geq & 2 \sqrt{\pi} \sqrt{\vert \conv (P) \vert}
 \nonumber \\
 & \geq & 2 \pi \sqrt{ T^2 - \frac{T^{4/3 + \epsilon'}}{\pi} } \nonumber \\
 & \geq & 2 \pi T -  \frac{\pi}{3} T^{\frac{1}{3} + \epsilon'}, \nonumber
\end{eqnarray}
where the first inequality uses the bound displayed in (\ref{cbv}). 
The second inequality is valid provided that $H(\epsilon')$
occurs and the third is true for all $T$ sufficiently high, since
$\epsilon' < 2/3$ and $3 < \pi$.

For $a \geq 0$, let $K_a$ denote the event that
$$
2^{a-1} T^{\frac{2}{3} + 2 \epsilon} \geq L_1 \geq 2^{a-2} T^{\frac{2}{3} + 2 \epsilon}.
$$
We claim that, for any $a \in \mathbb{N}$, 
\begin{equation}\label{piol}
 H \cap K_a  \subseteq \Big\{ \sum_{i=1}^{2m}{E_i^2} \geq 4 \pi^2 T +
 2^{2a-5} T^{\frac{2}{3} + 3 \epsilon} \Big\},
\end{equation}
where the collection of random variables $\big\{ E_i: i \in \{
1,\ldots, 2m \} \big\}$ was introduced in (\ref{rew}). To derive
(\ref{piol}), note that 
\begin{equation}\label{piola}
\sum_{i=1}^{m}{L_i^2} \geq L_1^2 + \frac{\Big(
\sum_{i=2}^{m}{L_i} \Big)^2}{m - 1}.
\end{equation}
From (\ref{piola}) and (\ref{rew}) follows 
\begin{equation}\label{piole}
\sum_{i=1}^{2m}{E_i^2} \geq \frac{m}{T}L_1^2 + \frac{\Big(
\sum_{i=2}^{m}{L_i} \Big)^2}{T}.
\end{equation}
Supposing the occurrence of the event $H (\epsilon') \cap K_a$, we
find that, for high enough values of $T$,
\begin{eqnarray}
\sum_{i=1}^{2m}{E_i^2} & \geq & \frac{T^{1/3 - \epsilon} - 1}{T}
2^{2a-4} T^{\frac{4}{3} + 4\epsilon} + \frac{1}{T} \Big( 2 \pi T -
\frac{\pi}{3} T^{\frac{1}{3} + \epsilon'}  - 2^{a-1} T^{\frac{2}{3} +
2 \epsilon} \Big)^2 \nonumber \\
 & \geq & 
2^{2a-4} T^{\frac{2}{3} + 3\epsilon} - 
2^{2a-4} T^{\frac{1}{3} + 4\epsilon}  
 +  4 \pi^2 T  - \frac{4\pi^2}{3} T^{\frac{1}{3} + \epsilon'}
 -  2^{a+1} \pi  T^{\frac{2}{3} + 2 \epsilon} \nonumber \\
 & \geq &  4 \pi^2 T   + 
2^{2a - 5} T^{\frac{2}{3} + 3\epsilon}. \nonumber
\end{eqnarray}
In the first inequality, (\ref{hcon}) and (\ref{piole}) were used,
while, in the third, the relations $\epsilon < 1/3$ and $\epsilon' - 3
\epsilon < 1/3$ were required.

From (\ref{piol}) and the definition of the random variables $E_i$ 
(with  $i \in \{ 1,\ldots, 2m \}$), it follows that
\begin{equation}\label{pionl}
\mathbb{P}(K_a \cap H) \leq \int_{l}^{\infty}{f_{2m-2}(y)dy},
\end{equation}
where $f$ denotes the density of the $\chi^2$-distribution specified
in (\ref{fchi}), and the lower limit of integration $l$ is equal to $4
\pi^2 T + 2^{2a - 5} T^{\frac{2}{3} + 3 \epsilon}$.
Using the fact that $m \leq T^{\frac{1}{3}}$ in estimating the
right-hand-side of (\ref{pionl}), we find that 
\begin{equation}\label{pionm}
\mathbb{P} \big( K_a \cap H(\epsilon') \big) \leq  \exp \Big\{ - 2 \pi^2 T  - 2^{2a-7}
T^{\frac{2}{3} + 3\epsilon} \Big\} .
\end{equation}
for all choices of $T$
exceeding a value that has no dependence on $a \geq 0$.
An extra factor of one-half multiplies the term of $
T^{\frac{2}{3} + 3\epsilon}$ in the last expression, in order to
compensate for the non-leading terms in the expression for the
$\chi^2$-density. By summing over $a \in \mathbb{N}$ in (\ref{pionm}), we deduce that,
for high values of $T$,
\begin{equation}\label{twnine}
\mathbb{P} \Big( \Big\{ L_1 \geq \frac{1}{4} T^{\frac{2}{3} + 2 \epsilon} \Big\}
\cap H(\epsilon') \Big) \leq  \exp \Big\{ - 2 \pi^2 T  - 2^{-8}
T^{\frac{2}{3} + 3 \epsilon } \Big\}.
\end{equation}
Allowing $I$ to denote the collection of those times $t \in [0,T]$ for
which 
\begin{displaymath}
 \Big\vert B \Big( t + T^{\frac{2}{3} + \epsilon} \Big) - B \big( t
 \big) \Big\vert \geq \frac{1}{2} T^{\frac{2}{3} + 2 \epsilon},
\end{displaymath}
note that the set $I$ is distributionally invariant under shifts of
$[0,T]$, and is non-empty if and only if the event $J(\epsilon)$
occurs. We denote by $I_1$ the event that there exists $t \in I$
satisfying $ \big\vert t - t' \big\vert \leq 1$, where the difference $t-t'$ is
being computed modulo $T$, as each other such will be. We find that
\begin{equation}\label{sfd} 
 \mathbb{P} \Big(  I_1 \cap \Big\{ \big\vert \enc
 (B) \big\vert \geq \pi T^2 \Big\} \Big\vert J (\epsilon) \cap \Big\{ \big\vert \enc
 (B) \big\vert \geq \pi T^2 \Big\} \Big) \geq \frac{2}{T}.
\end{equation}
If the event $I_1$ occurs, we have that
\begin{eqnarray}
 L_1 & = & \Big\vert B \Big( \frac{T}{m} + t' \Big) - B \big( t' \big)
 \Big\vert  \label{bjm} \\
  & \geq &  \Big\vert B \Big( t +  T^{\frac{2}{3} + \epsilon} \Big) - B \big( t \big)
 \Big\vert -  \Big\vert B \Big( \frac{T}{m} + t' \Big) -  B \Big( t +
 T^{\frac{2}{3} + \epsilon} \Big) \Big\vert -  \Big\vert B \big( t' \big) - B
 \big( t \big) \Big\vert \nonumber 
\end{eqnarray}
where $t \in I$ satisfies $\vert t - t' \vert \leq 1$.
Writing $f(T) = 2T^{\frac{1}{3} + 2\epsilon} + 1$ and $g(T) =
\frac{1}{8} T^{\frac{2}{3} + 2 \epsilon}$, note that for high values
of $T$, on the event
$Q_{f,g}$,
\begin{equation}\label{xlm}
  \Big\vert B \Big( \frac{T}{m} + t' \Big) -  B \Big( t +
 T^{\frac{2}{3}} \Big) \Big\vert \leq g \big( T \big).
\end{equation}
This is because the difference between the two arguments of $B$ in (\ref{xlm})
 satisfies
\begin{eqnarray}
  \Big\vert \Big( \frac{T}{m} + t' \Big) -  \Big( t +
 T^{\frac{2}{3}} \Big) \Big\vert & \leq & \bigg\vert \frac{T}{\lfloor
 T^{\frac{1}{3} - \epsilon} \rfloor }  - T^{\frac{2}{3} + \epsilon}
 \bigg\vert + \big\vert t' - t \big\vert \nonumber \\
 & \leq &  \bigg\vert \frac{T}{
 T^{\frac{1}{3} - \epsilon} - 1 }  - T^{\frac{2}{3} + \epsilon}
 \bigg\vert  + 1  \nonumber \\
 & \leq & T^{\frac{2}{3} + \epsilon} \sum_{j=1}^{\infty}{T^{-j \big(
 \frac{1}{3} - \epsilon \big)}} \, + \, 1 \leq f(T),
\end{eqnarray}
where the fact that $\epsilon \in (0,1/3)$ was used in the last inequality.
Note also that, on the event $Q_{1,g}$,
\begin{equation}\label{xln}
  \Big\vert B \big( t' \big) -  B \big( t \big) \Big\vert \leq g \big( T \big) .
\end{equation}
From (\ref{bjm}), (\ref{xlm}), (\ref{xln}) and the fact that $t \in I$, it follows that
\begin{displaymath}
 I_1 \cap Q_{f,g} \cap Q_{1,g} \subseteq \Big\{ L_1 \geq \frac{1}{4}
 T^{\frac{2}{3} + 2 \epsilon} \Big\}.
\end{displaymath}
Thus,
\begin{eqnarray}
 & &  \mathbb{P} \Big( \Big\{ L_1 \geq \frac{1}{4} T^{\frac{2}{3} + 2
 \epsilon} \Big\} \cap \Big\{ \big\vert \enc (B) \big\vert \geq \pi T^2
 \Big\} \Big) \label{sdx} \\
 & \geq &  \mathbb{P} \Big(  I_1 \cap \Big\{ \big\vert \enc (B)
 \big\vert \geq \pi T^2 \Big\} \Big) - \mathbb{P} \big( Q_{f,g}^c
 \big) - \mathbb{P} \big( Q_{1,g}^c \big) \nonumber \\ 
 & \geq & \frac{2}{T}  \mathbb{P} \Big( J \big( \epsilon \big) \cap \Big\{ \big\vert \enc (B)
 \big\vert \geq \pi T^2 \Big\} \Big) \nonumber \\
 & & \qquad - \, \, C T \exp \Big\{ - c T^{1 + 2
 \epsilon} \Big\} \, - \,
 C T \exp \Big\{ - \frac{1}{64 \cdot 128} T^{\frac{4}{3} + 4 \epsilon} \Big\}, \nonumber 
\end{eqnarray}
with $c \in (0,128^{-2})$.
In the second inequality, (\ref{sfd}) was used, as well as the
bound provided by Lemma \ref{remaftdef}. In applying this bound, we used the
fact that $f,g \geq 1$ for $T$ sufficiently large.
We find that
\begin{eqnarray}
 & & \mathbb{P} \Big( J \big( \epsilon \big) \cap \Big\{ \big\vert \enc (B)
 \big\vert \geq \pi T^2 \Big\} \Big) \nonumber \\
 & \leq & \frac{T}{2} \bigg[ \mathbb{P} \Big( \Big\{ L_1 \geq \frac{1}{4} T^{\frac{2}{3} + 2
 \epsilon} \Big\} \cap H \big( \epsilon' \big) \Big) + \mathbb{P}
 \Big( H (\epsilon')^c \cap \Big\{ \big\vert \enc(B) \big\vert \geq \pi
 T^2 \Big\} \Big)  \bigg] \nonumber \\
 & & \qquad \qquad + \, C T^2 \exp \Big\{ -c T^{1 + 2\epsilon} \Big\}
 \nonumber \\
 & \leq & \frac{T}{2} \exp \Big\{ - 2 \pi^2 T - 2^{-8} T^{\frac{2}{3}
 + 3 \epsilon} \Big\} \nonumber \\
 & & \qquad +  \, \frac{T}{2} \exp \Big\{ - 2 \pi^2 T -
 c_1 T^{\frac{1}{3} 
 +  2 ( \epsilon' - \epsilon)} \Big\} +  C T^2 \exp \Big\{-c T^{1 +
 2\epsilon} \Big\}
 \nonumber,  
\end{eqnarray}
where (\ref{sdx}) was used in the first inequality, the latter
requiring (\ref{twnine}) and an application of Lemma \ref{lemcl} (so
that we use the inequality $\epsilon' > 2 \epsilon$, which follows
from the assumption that $2 \epsilon' > 5 \epsilon + 1/3$.) 
From the fact that $5 \epsilon + 1/3 < 2 \epsilon'$ follows 
\begin{displaymath}
\mathbb{P} \Big( J \big( \epsilon \big) \cap \Big\{ \big\vert \enc (B)
 \big\vert \geq \pi T^2 \Big\} \Big) \leq  \exp \Big\{ - 2 \pi^2 T - 2^{-9} T^{\frac{2}{3}
 + 3 \epsilon} \Big\},
\end{displaymath}
for $T$ sufficiently high, as required. 

There remains the case where $\epsilon \in [1/5,1/3)$. Note that, for any
positive value for $\epsilon$, $J (\epsilon) \subseteq Q_{f,g}^c$,
with the choices 
\begin{displaymath}
  f(T) = T^{\frac{2}{3} + \epsilon}, \, g(T)= \frac{1}{4} T^{\frac{2}{3}
  + 2 \epsilon}
\end{displaymath}
being made. 
We find that, in this case,
\begin{eqnarray}
 & & \mathbb{P} \Big( J(\epsilon) \cap \Big\{  \big\vert \enc(B) \big\vert \geq \pi T^2
\Big\} \Big) \leq \mathbb{P} \big( J(\epsilon)) \leq \mathbb{P} \big(
Q_{f,g}^c \big) \label{gghh} \\ 
 & \leq & \frac{64}{\sqrt{\pi}} \Big( 2
 \sqrt{2} T^{-\frac{5\epsilon}{2}} + T^{-\frac{1}{6}- 2 \epsilon}  \Big)
  \exp \Big\{- \frac{1}{2048}
 T^{\frac{2}{3} + 3 \epsilon} \Big\}, \nonumber
\end{eqnarray}
where  Lemma \ref{remaftdef} was applied in the third
inequality. Given that $2/3 + 3 \epsilon > 1$ for a choice of
$\epsilon \in [1/5,1/3)$, we see that (\ref{gghh}) establishes  
the statement of the lemma for such values of $\epsilon$.
$\Box$ \\
{\bf Proof of Proposition \ref{ttwo}}
We claim that the following inclusion holds, for sufficiently high
values of $T$:
\begin{equation}\label{das}
 \bigcup_{\epsilon'} \Big\{ \exists  t \in [0,T]: \Big\vert B \Big( t + T^{\frac{2}{3} +
\epsilon'} \Big) - B \big( t \big) \Big\vert \geq T^{\frac{2}{3} + 2
\epsilon'} \Big\} \, \cap \,  Q_{\frac{1}{CT},1} \, \subseteq \, J \big( \epsilon \big),
\end{equation}
where the union on the left-hand-side is taken over values of
$\epsilon'$ satisfying $\vert \epsilon' - \epsilon \vert \leq
 \frac{T^{-5/3
 - \epsilon}}{2C \log T}$. 
To derive (\ref{das}), firstly set $V_{r} =  \frac{T^{-5/3 -
 r}}{2 C \log T}$ for $r>0$, and  note that, for sufficiently high values of
$T$, and for given $\epsilon > 0$ and any
$\epsilon'$ satisfying 
$$
\vert \epsilon' - \epsilon \vert \leq V_{\epsilon}
$$ 
we have that
\begin{equation}\label{dat}
 \Big\vert T^{\frac{2}{3} + \epsilon'} - T^{\frac{2}{3} + \epsilon}
 \Big\vert \leq \frac{1}{CT}.
\end{equation}
If the event on the left-hand-side of (\ref{das}) occurs, there exists
$\epsilon'$  satisfying $\vert \epsilon' - \epsilon \vert \leq T^{-5/3 -
\epsilon}/{\big( 2 C \log T \big)} $ and $t \in [0,T]$ for which 
\begin{equation}\label{daw}
\Big\vert B \Big( t + T^{\frac{2}{3} +
\epsilon'} \Big) - B \big( t \big) \Big\vert \geq T^{\frac{2}{3} + 2
\epsilon'}.
\end{equation}
Note that, provided that  the event on the left-hand-side of (\ref{das}) occurs,
\begin{eqnarray}
 & & \Big\vert B \Big( t + T^{\frac{2}{3} +
\epsilon} \Big) - B \big( t \big) \Big\vert \nonumber \\
 & \geq &  \Big\vert B \Big( t + T^{\frac{2}{3} +
\epsilon'} \Big) - B \big( t \big) \Big\vert -  \Big\vert B \Big( t + T^{\frac{2}{3} +
\epsilon'} \Big) - B \Big( t + T^{\frac{2}{3} + \epsilon} \Big)
\Big\vert \nonumber \\
 & \geq &   T^{\frac{2}{3} + 2 \epsilon'} - 1 
 =  T^{\frac{2}{3} + 2 \epsilon} + \big(   T^{\frac{2}{3} + 
 \epsilon'} - T^{\frac{2}{3} + \epsilon } \big) 
 \big( T^{\epsilon'} + T^{\epsilon} \big) - 1 
\geq \frac{1}{2}
 T^{\frac{2}{3} + 2 \epsilon}, \nonumber
\end{eqnarray}
where the final inequality, valid for high values of $T$, is due to
(\ref{dat}) and $\max \{\epsilon, \epsilon' \} < 1$.
In the second inequality, we have used (\ref{dat}), (\ref{daw}) and
the occurrence of $Q_{1/{CT},1}$. This establishes the inclusion (\ref{das}).

Note that
\begin{equation}\label{dau}
 \overline{J}_{\epsilon_0} \subseteq \bigcup_{i=1}^{N}{} \bigcup_{\epsilon' : \vert \epsilon' - \epsilon_i \vert \leq V_{\epsilon_i}} \Big\{ \exists  t \in [0,T]: \Big\vert B \Big( t + T^{\frac{2}{3} +
\epsilon'} \Big) - B \big( t \big) \Big\vert \geq T^{\frac{2}{3} + 2
\epsilon'} \Big\}, 
\end{equation}
where $\big\{ \epsilon_i : i \in \{1,\ldots, N \} \big\}$ denotes a
collection of values, each lying in $(\epsilon_0 , 1/3)$, for which 
\begin{displaymath}
 \Big[ \epsilon_0 , \frac{1}{3} \Big] \subseteq \bigcup_{i=1}^{N}{ \Big(
 \epsilon_i - V_{\epsilon_i}, \epsilon_i + V_{\epsilon_i}   \Big) }.
\end{displaymath}
As such, we may choose 
\begin{equation}\label{day}
 N \leq C T^2 \log T .
\end{equation}
From (\ref{das}) and (\ref{dau}) follows
\begin{displaymath}
 \overline{J}_{\epsilon_0} \cap \Big\{ \big\vert \enc (B) \big\vert
 \geq \pi T^2 \Big\} \subseteq \Big( \bigcup_{i=1}^{N}{J \big( \epsilon_i
 \big)} \cap  \big\{ \big\vert \enc (B) \big\vert
 \geq \pi T^2 \big\} \Big) \cup Q_{\frac{1}{CT},1}^c.
\end{displaymath}
From Lemma \ref{tone}, the bound on $\mathbb{P}\big( Q_{f,g} \big)$
given by Lemma \ref{remaftdef}, and (\ref{day}), it follows that
\begin{eqnarray}
 & & \mathbb{P} \Big( \overline{J}_{\epsilon_0} \cap \Big\{ \big\vert
 enc(B) \big\vert \geq \pi T^2 \Big\} \Big) \nonumber \\
 & \leq &
  C T^2 \log T \exp{ \Big\{ - 2\pi^2 T - 2^{-9} T^{\frac{2}{3} + 3
 \epsilon_0 } \Big\} } \nonumber \\
 & &   + \,  \Big( \frac{32\sqrt{2}}{\sqrt{\pi}} \sqrt{C} T^{\frac{3}{2}} +
 \frac{16}{\sqrt{\pi}} T^{\frac{1}{2}} \Big)  \exp \Big\{-\frac{C}{128} T \Big\}. \nonumber
\end{eqnarray}
The second term here is negligible provided that $C > 20 \pi^2$. 
Choosing $C_3 = 2 C$ yields the statement in the Proposition. $\Box$
\end{subsection}
\end{section}
\begin{section}{Proof of Theorem \ref{tgg}}\label{secthr}
We require the following result.
\begin{lemma}\label{lone}
Let $K$ denote a planar compact convex set whose area exceeds
$\pi T^2$. Then
\begin{displaymath}
   \arcl \big( \partial K \big)^2 > 4 \pi^2 T^2 + \pi^2 \big( \rout
   ( K ) - \rin(K)\big)^2.
\end{displaymath}
\end{lemma}
{\bf Proof} The result is implied by Bonnesen's inequality, as it is
stated in the Theorem of Subsection $1.3.1$, on page $3$ of \cite{burago}.
 $\Box$ \\
{\bf Proof of Theorem \ref{tgg}}. 
We will prove that, for $\epsilon \in (0,1/6)$ and all sufficiently
high $T$,
\begin{eqnarray}
 & & \mathbb{P} \Big( \Big\{  \rout ( B )  -  \rin ( B ) >  T^{\frac{2}{3} + \epsilon} \Big\} \cap 
\Big\{ \big\vert \enc(B) \big\vert \geq \pi T^2 \Big\} \Big)
\label{xgh} \\
 & \leq &  
 \exp \Big\{ - 2 \pi^2 T -  c_4 T^{\frac{1}{3} + 2 \epsilon} \Big\}, \nonumber
\end{eqnarray}
for any $c_4 \in (0,\pi^2/{32})$. In doing so,
we will make use of the polygon $P$, in
the case where $m = \lfloor T^{\frac{1}{3}} \rfloor$ and $t' = 0$.  
\begin{definition}
Let $\epsilon_0,\epsilon_1,\epsilon_2,\epsilon_3,\epsilon_4$  be
positive constants that satisfy the following bounds:
$\epsilon_0 < 1/3$, $2
\epsilon_2 > \epsilon_1$, $\epsilon_1 < 1/6$, $\epsilon_2 < 1/3$,
$\epsilon_3 > {\epsilon_2}/4$, $\epsilon_4 > 2 \epsilon_0$.
\begin{itemize}
\item Let $H_1$ denote the event that 
$$
 66 \pi^3 T^2 \geq \big\vert \conv(P) \big\vert \geq \pi T^2 -
T^{\frac{4}{3} + \epsilon_1}.
$$
\item Let $H_2$ be given by
$$
H_2 = \Big\{ \rout \big(  \conv (P) \big) - \rin \big(  \conv (P) \big)  \leq T^{\frac{2}{3} +
\epsilon_2} \Big\}.
$$
\item Let $H_3$ be given by
$$
H_3 = \Big\{ \sup_{k \in \{ 0,\ldots, m-1 \} }d \Big( \partial \big(
conv P \big) , B \Big( \frac{kT}{m} \Big) \Big) \leq T^{\frac{2}{3} + \epsilon_3} \Big\}.
$$
\item Let $H_4$ denote the event that, for each $i \in \{ 0, \ldots, m-1
\}$,
$$
\sup_{t \in [0, T/m]}{\Big\vert B \Big( \frac{(i-1)T}{m} + t \Big) - B
\Big( \frac{(i-1)T}{m} \Big) \Big\vert \leq T^{\frac{2}{3} + \epsilon_4}}.
$$ 
\end{itemize}
\end{definition}
Note that
\begin{eqnarray}
 & & \mathbb{P} \Big( H_1^c \cap \Big\{ \big\vert \enc(B) \big\vert  
\geq \pi T^2 \Big\} \Big) \label{sone} \\
 & \leq &  \mathbb{P} \Big(  \Big\{ \big\vert \conv(P) \big\vert < \pi T^2 -
T^{\frac{4}{3} + \epsilon_1} \Big\} \cap \Big\{ \big\vert \enc(B) \big\vert  
\geq \pi T^2 \Big\} \Big) \nonumber \\
 & & \qquad  + \,  \mathbb{P} \Big( \big\vert 
\conv(P) \big\vert > 66 \pi^3 T^2 \Big). \nonumber
\end{eqnarray}
The first event after the inequality in (\ref{sone}) coincides with
$H(\epsilon_1)^c \cap  \big\{ \big\vert \enc(B) \big\vert  
\geq \pi T^2 \big\}$ as it appears in Lemma \ref{lemcl}, with $t'$ set
equal to zero. Applying this lemma with the choice $\epsilon = 0$, 
we deduce that, for sufficiently high values of $T$, the first term
on the right-hand-side of (\ref{sone}) is bounded above by
\begin{equation}\label{ston}
 \exp \Big\{- 2 \pi^2 T - c_1 T^{\frac{1}{3} + 2\epsilon_1} \Big\}.
\end{equation}
Note also that, if $ \big\vert 
\conv( P) \big\vert > 66 \pi^3 T^2$, then
\begin{equation}\label{nlu}
 \sup_{s,t \in [0,T]}{\big\vert B(s) - B(t) \big\vert^2} \geq {\rm diam}
 \big( \conv(P) \big)^2 \geq \frac{4}{\pi} \big\vert \conv(P)
 \big\vert \geq 4 \cdot 66 \pi^2 T^2,
\end{equation}
the first inequality being valid because the endpoints of the longest
diameter of $\conv(P)$ are vertices $B(t_1)$ and $B(t_2)$ of $P$, the
second being the standard isoperimetric inequality. Note that
(\ref{nlu}) implies the occurrence of $Q_{T,T \pi \sqrt{8 \cdot
33}}^c$. 
Thus, we learn from Lemma \ref{remaftdef} that
\begin{equation}\label{sten}
 \mathbb{P} \Big( \big\vert 
\conv( P) \big\vert > 6 \pi^3 T^2 \Big) \leq \mathbb{P} \Big(
Q_{T,T \pi \sqrt{8 \cdot 33}}^c \Big) \leq C \exp \Big\{ - \frac{33}{16} \pi^2 T \Big\}, 
\end{equation}
By
(\ref{sone}), the bound (\ref{ston}) on the first term on its
right-hand-side, and (\ref{sten}), we have that, for high values of
$T$,
\begin{equation}\label{sgt}
 \mathbb{P} \Big( H_1^c \cap \Big\{ \big\vert \enc(B) \big\vert  
\geq \pi T^2 \Big\} \Big) \leq \exp \Big\{- 2 \pi^2 T - \frac{c_1}{2} T^{\frac{1}{3} + 2\epsilon_1} \Big\}.
\end{equation}
We assert that
\begin{equation}\label{stwo}
\mathbb{P} \Big( H_1 \cap H_2^c \Big)
\leq \exp \Big\{ - 2 \pi^2 T - c_2 T^{\frac{1}{3} + 2 \epsilon_2} \Big\},
\end{equation}
for any choice of $c_2 \in \big( 0,\pi^2/2 \big)$. 
To show this, we choose $\delta$ to satisfy $\delta < 2/3$ and $\delta
\in (\epsilon_1,2\epsilon_2)$.
Note firstly that,
since $\delta \in (\epsilon_1, 2/3)$, 
$$
\pi T^2 - T^{\frac{4}{3} + \epsilon_1} \geq \pi \big( T -
T^{\frac{1}{3} + \delta} \big)^2,
$$
for $T$ sufficiently high.
The occurrence of $H_1$ therefore implies that
\begin{equation}\label{qsa}
 \big\vert  \conv( P )  \big\vert  \geq   \pi \big( T - T^{\frac{1}{3}
 + \delta} \big)^2 \, ,
\end{equation}
whereas, on the event $H_2^c$,
\begin{displaymath}
 \rout \big( \conv(P) \big) -  \rin \big( \conv(P) \big) >  T^{\frac{2}{3} + \epsilon_2}.
\end{displaymath}
Using the fact that $\delta < 2/3$,
we may apply Lemma \ref{lone} to find a lower bound on the arclength
of the convex hull of the polygon $P$ in this eventuality:
\begin{equation}\label{jpo}
\arcl \big( \conv( P ) \big)^2 \geq  4 \pi^2 \Big(  T - T^{\frac{1}{3}
+ \delta} \Big)^2  + \pi^2 T^{\frac{4}{3} + 2 \epsilon_2}. 
\end{equation}
From (\ref{jpo}) follows
\begin{eqnarray}
 \arcl \big( \conv(P) \big) & \geq & 2 \pi  \big(  T - T^{\frac{1}{3}
+ \delta} \big) \sqrt{1 + \frac{1}{4} T^{\frac{-2}{3} + 2 \epsilon_2}}
\nonumber \\
 & \geq &  2 \pi  \Big(  T - T^{\frac{1}{3}
+ \delta} \Big)  \bigg(  1 +  \Big( \frac{1}{8} - o(1) \Big) T^{\frac{-2}{3}
+ 2 \epsilon_2 } \bigg), \label{wfd}  
\end{eqnarray}
where the facts that $\delta < 2/3$ and  $\epsilon_2 < 1/3$  were used
in successive inequalities. Since $2 \epsilon_2 > \delta$, we deduce that
\begin{equation}\label{qq}
\arcl \big( \conv( P ) \big)
  \geq  2 \pi  T + c 
 T^{\frac{1}{3} + 2 \epsilon_2},
\end{equation}
for any $c \in (0,\pi/4)$.
Recalling Definition \ref{defnsix}, we find that
\begin{eqnarray}
 \sum_{i=1}^{m}{L_i^2} & \geq & \frac{\Big( \sum_{i=1}^{m}{L_i}
 \Big)^2}{m} \geq  \frac{\Big( 
\arcl \big( \conv( P ) \big) \Big)^2}{m} \nonumber \\
 & \geq &  \frac{\Big( 2 \pi T + c T^{\frac{1}{3} + 2 \epsilon_2}
 \Big)^2}{m} \geq  \frac{ 4 \pi^2 T^2 + 4 \pi c T^{\frac{4}{3} + 2 \epsilon_2}
 }{m}. \nonumber
\end{eqnarray}
Thus,
\begin{equation}\label{anotime}
\sum_{i=1}^{2m}{E_i^2} \geq  4 \pi^2 T + 4 \pi c T^{\frac{1}{3} + 2
\epsilon_2 },
\end{equation}
where the collection of random variables $\big\{ E_i : i \in \{1,\ldots,2m
\} \big\}$ was introduced in (\ref{rew}).
The left-hand-side of (\ref{anotime}) having the $\chi^2$-distribution
with $2m - 2$ degrees of freedom, and $m$ being at most $T^{1/3}$, it
follows by (\ref{fchi}) that the probability of the occurrence of (\ref{anotime}) is at most 
$$
 (CT)^{\frac{1}{2}T^{1/3}} \exp \Big\{ - 2 \pi^2 T - 2 \pi
 c T^{1/3 + 2\epsilon_2}  \Big\} \, + \exp\big\{-CT/4 \big\},
$$
for any constant $C$.
Since $c$ may be chosen to lie arbitrarily close to $\pi/4$, we have shown that (\ref{stwo}) holds.

We will now show that
\begin{equation}\label{qsb}
 \mathbb{P} \Big( H_1 \cap H_2 \cap H_3^c \Big) \leq \exp \Big\{ - 2
 \pi^2 T - c_3 T^{\frac{1}{2}} \Big\},
\end{equation}
for any $c_3 \in \big(0,2\pi \big)$. We do so by proving the inclusion
\begin{equation}\label{fifsiha}
 H_1 \cap H_2 \cap H_3^c \subseteq \Big\{ \sum_{i=1}^{m}{L_i} - \arcl
 \big( \partial (\conv P) \big) > T^{\frac{1}{2}} \Big\}, 
\end{equation}
for which purpose,
we require a lemma. 
\begin{lemma}\label{apl}
Let $P$ denote a planar polygon, with vertex set $\big\{ p_i : i \in
\{1,\ldots,m \} \big\}$ (so that its arclength $\arcl(P)$ is given by
the sum of its edge-lengths $\vert p_{i+1} - p_i \vert$ for $i \in \{1,\ldots,m\}$). 
Then
\begin{displaymath}
 \arcl \big( P \big) \geq \arcl \big( \partial (\conv P) \big) \, + \,  \big(
 \sqrt{5} - 2 \big) \min
 \Big\{ \frac{2R^2}{Q}, R \Big\},
\end{displaymath}
where $R = \sup_{i \in \{1,\ldots,m\}} d\big( p_i, \partial ( \conv P
) \big)$, and where $Q$ is equal to
the supremum of the lengths of line segments in $\partial \big( \conv
P \big)$.
\end{lemma}
\noindent{\bf Proof.} 
Let $\rm L$ denote the collection of line segments that
comprise $\partial \big( \conv P \big)$. 
Write $\partial_{\ext}(P)$ for the exterior boundary of the polygon
$P$. That is, $\partial_{\ext}(P) = \partial \big( \enc(P)^c
\big)$. Note that the set $\conv(P) \setminus \enc(P)$ is comprised of
a finite number of connected components. The boundary of each
component consists of the union of a line sement $l \in \rm L$, and a
polygonal path lying in $\partial_{\ext}(P)$ whose endpoints coincide
with those of $l$. We denote this polygonal path by $l_{\ext}$. In the
case that $l \in \rm L$ does not arise from any such component, we set
$l_{\ext} = l$. Note that, for any pair $l',\tilde{l} \in \rm L$, the
set ${l'}_{\ext} \cap \tilde{l}_{\ext}$ has at most finitely many
elements. Note further that $\vert l_{\ext} \vert \geq \vert l \vert$
for each $l \in \rm L$, where we use $\vert \cdot \vert$ to denote the
length of a line segment, or of a finite union of line segments.

Let $k \in \{1,\ldots,m\}$
satisfy 
\begin{equation}\label{rope}
 R = d \Big( p_k, \partial \big( \conv P \big) \Big).
\end{equation}
We distinguish two cases, according to whether or not 
\begin{equation}\label{tcs}
 d \Big( p_k , \partial_{\rm ext} \big( P \big) \Big)
 \Big) \leq
 \frac{R}{2}.
\end{equation}
Supposing that (\ref{tcs}) holds, let $x \in \partial_{\ext}\big(
P \big)$ satisfy 
\begin{equation}\label{yyu}
  d \Big( x , p_k \Big)
 \Big) \leq
 \frac{R}{2} .
\end{equation}
Let 
$l^* = \big[ p_i , p_j \big]$
denote the element of $\rm L$ for which $x \in l^*_{\ext}$.
By using (\ref{rope}) and (\ref{yyu}), we find that
\begin{equation}\label{bas}
  d \big( x , l^* \big) >  \frac{R}{2} .
\end{equation} 
Note also that
\begin{displaymath}
  \big\vert l^*_{\ext} \big\vert  \geq 
\big\vert x -  p_i
\big\vert + 
\big\vert p_j - x \big\vert 
 \geq  
\big\vert q -  p_i
\big\vert + 
\big\vert q -  p_j \big\vert, 
\end{displaymath}
where $q$ denotes the point at distance $\frac{R}{2}$ from $l^*$
whose projection onto this line segment is its midpoint. In the latter
inequality, we invoked (\ref{bas}). Thus,
\begin{eqnarray}
& & \big\vert l^*_{\ext} \big\vert \geq  \sqrt{\big\vert 
 p_j -  p_i \big\vert^2 + R^2 } \label{fifeha} \\ 
 & = &  \big\vert p_j - p_i \big\vert \sqrt{1 +
\frac{R^2}{\big\vert p_j - p_i \big\vert^2 }}. \label{srq}
\end{eqnarray}
If 
\begin{equation}\label{qse}
 \big\vert p_j -  p_i
 \big\vert > 2 R,
\end{equation}
then 
\begin{equation}\label{gopg} 
   \big\vert l^*_{\ext} \big\vert  \geq  \big\vert p_i - p_j
   \big\vert 
 + 2 \big( \sqrt{5} - 2 \big) \frac{R^2}{ \big\vert p_j - p_i \big\vert} 
   \geq  \big\vert p_i - p_j \big\vert +  2 \big( \sqrt{5} - 2 \big) \frac{R^2}{Q} . 
\end{equation}
In the first inequality of (\ref{gopg}), we used 
(\ref{qse}) and the fact that 
\begin{displaymath}
\sqrt{1 + x} \geq 1 +
 2 \big( \sqrt{5} - 2 \big) x
\end{displaymath} 
for $x \in [0,1/4]$ to bound below the term
appearing in (\ref{srq}). In the second inequality of  (\ref{gopg}), we used the
fact that $\vert l^* \vert \leq Q$ (which follows from the definition
of $Q$). 
Thus, provided that (\ref{tcs}) and (\ref{qse}) hold, we have that 
\begin{eqnarray}
 \arcl \big( P \big) & \geq & \sum_{l \in \rm L}{\vert l_{\ext} \vert} \geq \sum_{l \in
 {\rm L} \setminus \{ l^* \} }{\vert l \vert} \, + \, \big\vert p_i
  - p_j \big\vert \, + \,
  2\big( \sqrt{5} - 2 \big) \frac{R^2}{Q} \label{gopl} \\
 & = & \arcl \big( \partial ( \conv P ) \big) \, + \,    2\big( \sqrt{5} - 2 \big)
 \frac{R^2}{Q}. \nonumber 
\end{eqnarray}
where the second inequality is 
due to the fact that $ \vert l_{\ext} \vert \geq \vert l
\vert$ for each $l \in L$.
Note that, if (\ref{qse}) fails, then
\begin{equation}\label{sixtw}
   \big\vert l^*_{\ext} \big\vert    - \big\vert p_i - p_j \big\vert 
 \geq   \sqrt{\big\vert p_i - p_j \big\vert^2 +
 R^2 }  - \big\vert p_i - p_j  \big\vert  
  \geq  \big( \sqrt{5} - 2 \big) R, 
\end{equation}
where the first inequality follows from (\ref{fifeha}) and where the
second is obtained by minimising the value of the term  
$\sqrt{\big\vert p_i - p_j \big\vert^2 +
 R^2 }  - \big\vert p_i - p_j  \big\vert$  
subject
to the constraint that (\ref{qse}) does not hold. The inequality
arising from
(\ref{sixtw}) may replace (\ref{gopg}) in deriving (\ref{gopl}). We
have established the statement of the lemma in the event that
(\ref{tcs}) holds. 
In the other case, there exists a union $\hat{l}$ of line segments in
$P$ such that $p_k \in \hat{l}$, $\vert \hat{l} \vert
\geq \frac{R}{2}$, and $\hat{l} \cap
\partial_{\ext} \big( P \big) = \emptyset$. Given that $\big\vert
\partial_{\ext}\big( P \big) \big\vert \geq \arcl \big( \partial ( \conv P
 ) \big)$ (which follows from $\vert l_{\ext} \vert \geq \vert l \vert$ for each
$l \in \rm L$), we find that
\begin{displaymath}
 \arcl \big( P \big) \geq  \big\vert \partial_{\ext}\big( P \big)
 \big\vert \, + \, \big\vert \hat{l} \big\vert \geq \arcl\big(
 \partial (\conv P)
 \big) \, + \, \frac{R}{2} .
\end{displaymath}
Given that $\frac{1}{2} > \sqrt{5} - 2$, this 
establishes the statement of the lemma in the case that
(\ref{tcs}) fails, thereby completing the proof of the lemma. $\Box$

In applying Lemma \ref{apl}, we firstly find an upper bound on the
quantity $Q$ appearing in its statement, for the polygon $P$ under discussion.
To be specific, we now show that, provided that $H_1 \cap H_2$
occurs, then, for any line segment $L$ in the convex boundary of $P$,
\begin{equation}\label{qsc}
 \big\vert L \big\vert \leq C_4 T^{\frac{5}{6} + \frac{\epsilon_2}{2}},
\end{equation}
where $C_4$ is any constant exceeding $2^{3/2} (66)^{1/4} \sqrt{\pi}$.
Indeed, if $H_1 \cap H_2$ occurs, there exist two planar circles $J_1$
and $J_2$, satisfying
\begin{displaymath}
\rad(J_1)  \in \Big[ \frac{T}{2}, \pi \sqrt{66} T \Big], \, \rad (J_2) \in \Big[
\rad(J_1),\rad(J_1) + T^{\frac{2}{3} + \epsilon_2} \Big], 
\end{displaymath}
and 
$\enc \big( J_1 \big) \subseteq \enc \big( J_2 \big)$,
such that $\partial \big( \conv P \big)$
lies in $\enc(J_2) \setminus \enc(J_1)$. 

The length of the longest line segment $L$ that lies in  $\enc(J_2) \setminus
\enc(J_1)$ for a pair of circles $J_1$ and $J_2$ satisfying these
conditions is attained when $\rad(J_1)$ and $\rad(J_2)$ are maximal,
with the relative positions of $J_1,J_2$ and $L$ resembling the
picture in Figure 6. Thus, this length is bounded above
by
\begin{eqnarray}
 & & 2 \sqrt{\Big( \pi \sqrt{66} T + T^{\frac{2}{3} + \epsilon_2 } \Big)^2
 - 66 \pi^2 T^2 } \nonumber \\
 & \leq & 2 \sqrt{2 \pi \sqrt{66} T^{\frac{5}{3} + \epsilon_2} +
 T^{\frac{4}{3} + 2 \epsilon_2} } \leq  2^{3/2} (66)^{1/4} \sqrt{\pi}
 \big( 1 + o(1) \big) T^{\frac{5}{6} + \frac{\epsilon_2}{2}}, 
\end{eqnarray} 
the second inequality requiring $\epsilon_2 < 1/3$.
This establishes (\ref{qsc}).

To apply Lemma \ref{apl} and derive (\ref{fifsiha}), note 
that the occurrence of $H_3^c$ implies that $R > T^{\frac{2}{3}+
\epsilon_3}$. It follows from (\ref{qsc}) that, in the event of $H_1
\cap H_2$, we have the inequality, $Q \leq C_4 T^{\frac{5}{6} +
\frac{\epsilon_2}{2}}$. Making use of the assumption that $\epsilon_3
> \epsilon_2/4$, we indeed obtain (\ref{fifsiha}) from Lemma \ref{apl}.

Suppose that the quantity $\delta$ chosen before (\ref{qsa}) satisfies
$\delta < 1/6$ (at this point, we require that $\epsilon_1 < 1/6$). Applying the standard isoperimetic inequality and
using (\ref{qsa}), which holds provided that $H_1$ occurs, we find that, on $H_1$,
$$
\arcl \big( \partial ( \conv P ) \big) \geq 2\pi T - 2 \pi T^{\frac{1}{3} + \delta}.
$$ 
From (\ref{gopl}), whose validity we have established whether or not
(\ref{qse}) holds, we deduce that
$$
\sum_{i=1}^{m}{L_i} \geq 2 \pi T + \big( 1 - o(1) \big) T^{\frac{1}{2}},
$$
since $\delta < 1/6$.
By a reprise of the argument that follows (\ref{qq}),
the probability of this last
inequality is at most $\exp \Big\{ - 2 \pi^2 T -
 2\pi \big( 1 - o(1) \big) T^{\frac{1}{2}}\Big\}$. We have now  established (\ref{qsb}). 

Given that $\epsilon_0 < \min\big\{ 1/3, \epsilon_4/2 \big\}$, Lemma \ref{smalllemma} implies that
$\overline{J}_{\epsilon_0}^c \subseteq H_4$. 
Let $\psi(T) = T^{\frac{2}{3}} \big( T^{\epsilon_2} + 2 T^{\epsilon_3}
+ 2 T^{\epsilon_4} \big)$.
It is straightforward that
\begin{equation}\label{qsg}
H_2 \cap H_3 \cap H_4 \subseteq \Big\{ \rout ( B ) -  \rin ( B ) \leq
\psi(T) \Big\}.
\end{equation}
From (\ref{qsg}), it follows that 
\begin{eqnarray}
& & \mathbb{P} \Big( \Big\{ \rout ( B )  -  \rin ( B ) > \psi(T) \Big\} \cap 
\Big\{ \big\vert \enc(B) \big\vert \geq \pi T^2 \Big\} \Big) \nonumber
\\
& \leq & \mathbb{P} \Big( ( H_2 \cap H_3 \cap H_4 )^c \cap 
\Big\{ \big\vert \enc(B) \big\vert \geq \pi T^2 \Big\} \Big) \nonumber
\\
& \leq & \mathbb{P} \Big( H_1^c  \cap 
\Big\{ \big\vert \enc(B) \big\vert \geq \pi T^2 \Big\} \Big) +
\mathbb{P} \big( H_1 \cap H_2^c \big) \nonumber \\
 & &  + \, \mathbb{P} \big( H_1 \cap H_2
\cap H_3^c \big)
+  \mathbb{P} \Big( \overline{J}_{\epsilon_0} \cap \Big\{ \big\vert \enc(B) \big\vert \geq
\pi T^2 \Big\}  \Big) \nonumber \\
& \leq &  \exp \Big\{ - 2 \pi^2 T - \frac{c_1}{2} T^{\frac{1}{3} + 2 \epsilon_1} \Big\} +
\exp \Big\{ - 2 \pi^2 T - c_2 T^{\frac{1}{3} + 2 \epsilon_2} \Big\}  
\label{eighthal} \\
 & & + \exp \Big\{ - 2 \pi^2 T - c_3 T^{\frac{1}{2}} \Big\} \, + C_3 T^2 \log T \exp \Big\{ - 2 \pi^2 T -
 2^{-9} T^{\frac{2}{3} + 3 \epsilon_0} \Big\}, \nonumber
\end{eqnarray}
where, in the second inequality, we used $H_4^c \subseteq
\overline{J}_{\epsilon_0}$. The bounds (\ref{sgt}), (\ref{stwo}),
(\ref{qsb}) and that provided by Proposition \ref{ttwo} were used for
the successive terms to obtain the third inequality. 
For $\epsilon \in (0,1/6)$, we set  $\epsilon_1 = \epsilon + \hat{\epsilon}, \epsilon_2 =
\epsilon_4 = \epsilon, \epsilon_3 = \epsilon/4 + \hat{\epsilon}$ and
$\epsilon_0 = \epsilon/2 - \hat{\epsilon}$, where $\hat{\epsilon}>0$
is chosen small enough that the restrictions on the parameters $\big\{
\epsilon_i : i \in \{0,\ldots,4\} \big\}$ are satisfied.
We obtain
\begin{eqnarray}
 & & 
\mathbb{P} \Big( \Big\{  \rout ( B )  -  \rin ( B ) > 3 T^{\frac{2}{3} + \epsilon} \Big\} \cap 
\Big\{ \big\vert \enc(B) \big\vert \geq \pi T^2 \Big\} \Big) \nonumber
\\
 & \leq & 
 \exp \Big\{ - 2 \pi^2 T - \big( c_2 - o(1) \big) T^{\frac{1}{3} + 2
 \epsilon} \Big\}, \nonumber
\end{eqnarray}
since $\epsilon \in (0,1/6)$. From this follows (\ref{xgh}), the value
of $c_4$ in (\ref{xgh}) arising because we may choose $c_2$ in (\ref{stwo})
to be  slightly
less than $\pi^2/2$.
Note that
\begin{eqnarray}
 & &  \Big\{ \rin(B) < T - T^{\frac{2}{3} + \epsilon} \big\} \cap \big\{
 \vert \enc(B) \vert \geq \pi T^2 \Big\} \label{pffp} \\
 & \subseteq & \Big\{ \rout(B) -
 \rin(B) > T^{\frac{2}{3} + \epsilon} \Big\}, \nonumber
\end{eqnarray}
because $\rout(B) \geq T$ if $\vert \enc(B) \vert \geq \pi T^2$.

We obtain
\begin{eqnarray}
 & & \mathbb{P} \Big( \rin \big( X \big) < T - T^{\frac{2}{3} +
 \epsilon} \Big) \nonumber \\
 & = & \frac{\mathbb{P} \Big( \Big\{ \rin \big( B \big) < T - T^{\frac{2}{3} +
 \epsilon} \Big\} \cap \Big\{ \big\vert \enc(B)\big\vert \geq \pi T^2 \Big\}\Big)}{\mathbb{P} \Big( \big\vert \enc(B) \big\vert \geq \pi
 T^2 \Big)} \nonumber \\
 & \leq &  \exp \big\{ - c_4 T^{\frac{1}{3} + 2 \epsilon} + C_1
 T^{\frac{1}{3}} \log T \Big\}, \nonumber
\end{eqnarray}
the inequality by means of (\ref{xgh}), (\ref{pffp}) 
and Lemma \ref{ltwo}. We have proved the first part of the
theorem. 

To derive the second statement, note that
\begin{eqnarray}
 & & \Big\{ \rout(B) > T + 2 T^{\frac{2}{3} + \epsilon} \Big\} \cap \Big\{
 \vert \enc(B) \vert \geq \pi T^2 \Big\} \nonumber \\
 & \subseteq &  \Big\{ \rin (B) > T + T^{\frac{2}{3} + \epsilon}
 \Big\}  \label{jjk} \\ 
 & & \cup \, \Big\{ \rout(B) - \rin (B) > T^{\frac{2}{3} + \epsilon} \Big\} \cap \Big\{
 \vert \enc(B) \vert \geq \pi T^2 \Big\}. \label{ppl}
\end{eqnarray}
Noting that the inequality $\rin (B) > T + T^{\frac{2}{3} + \epsilon}$
implies that $\big\vert \enc(B) \big\vert \geq \pi T^2 + 2\pi
T^{\frac{5}{3} + \epsilon}$, we may use Proposition \ref{propexcarea} to bound
the probability of the event in (\ref{jjk}). We obtain that, for any $\alpha \in (0,1/3)$,
\begin{equation}\label{plk}
 \mathbb{P} \Big(  \rin (B) > T + T^{\frac{2}{3} + \epsilon} \Big)
 \leq \exp \Big\{ - 2 \pi^2 T - T^{\frac{1}{3} + \alpha} \Big\},
\end{equation}
provided that $T$ is chosen to be sufficiently  high.

Using (\ref{xgh}) to bound the probability of the event in 
(\ref{ppl}), as well as (\ref{plk})
we find that, for $\epsilon \in (0,1/6)$
and for all sufficiently high values of $T$,
\begin{eqnarray}
 & & \mathbb{P} \Big( \Big\{ \rout(B) > T + 2 T^{\frac{2}{3} + \epsilon} \Big\} \cap \Big\{
 \vert \enc(B) \vert \geq \pi T^2 \Big\} \Big) \label{hhk} \\ 
 & \leq &   \exp \Big\{ - 2 \pi^2 T - \frac{c_4}{2} T^{\frac{1}{3} + 2 \epsilon}
  \Big\}. \nonumber 
\end{eqnarray}
From (\ref{hhk}) and Lemma \ref{ltwo}, we deduce the second statement of the theorem. $\Box$ \\
\begin{lemma}\label{lemfi}
Let $K$ denote a planar compact set. In the case that $\rin (K) >
\frac{1}{2} \rout(K)$, we have that
\begin{displaymath} 
 \mathcal{L} \big( \conv K \big) \leq 4 \sqrt{\rin \big( K \big)  \big( \rout ( K )  -
 \rin ( K ) \big) }. 
\end{displaymath}
\end{lemma}
{\bf Proof} Given that $K$ is compact, we may locate circles $J_1$ and $J_2$
satisfying 
\begin{eqnarray}
 & & \enc \big( J_1 \big) \subseteq K \subseteq \enc \big( J_2 \big) , \, \, \rad
\big( J_2 \big) = \rout \big( K \big)  \nonumber \\
 & & \qquad \qquad \textrm{and} \, \, \,  \rad
\big( J_1 \big) = \rin \big( K \big).\label{denone}
\end{eqnarray}
Note that any line segment $L$ lying in $\partial \big( \conv K \big)$ satisfies
\begin{equation}\label{dentwo}
L \subseteq  \intn \big( \enc \big( J_1 \big) \big)^c \cap \enc \big( J_2 \big),
\end{equation}
where $\intn (A)$ denotes the interior of the set $A$. 
Indeed, if $\intn \big( \enc ( J_1 ) \big) \cap L \not= \emptyset$, there exists
a point lying in $\enc
\big( J_1 \big)$ and in the half-plane whose boundary contains $L$ and
is disjoint from $\conv \big( K \big)$. This point does not lie in $K$, 
implying that $\enc \big( J_1 \big) \not\subseteq K$. The
endpoints of $L$ lie in $\enc \big( J_2 \big)$, and thus, so does $L$,
by the convexity of $\enc \big( J_2 \big)$. 
\begin{figure}\label{dfour}
\begin{center}
\includegraphics[width=0.5\textwidth]{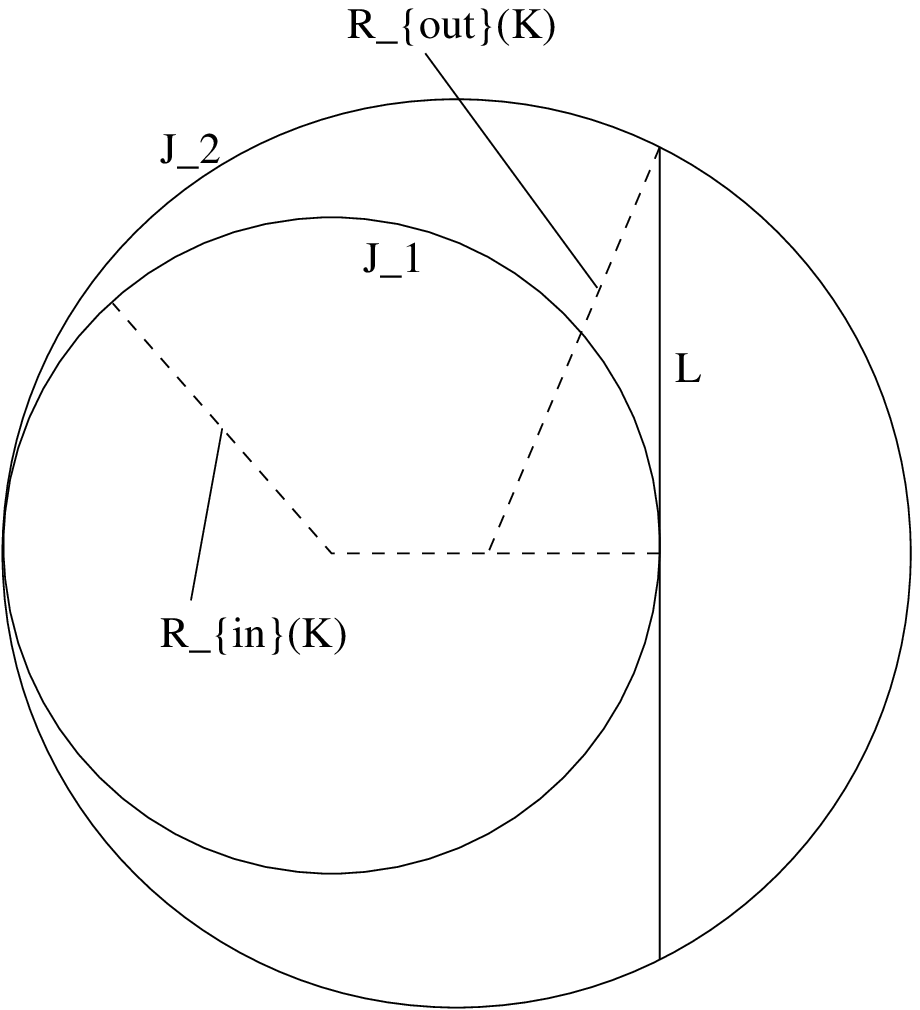}
\end{center}
\caption{Computing the length $\big\vert L \big\vert$ in (\ref{lenell})}.
\end{figure}

Given that  $\rin (K) >
\frac{1}{2} \rout(K)$, the supremum of the lengths of line segments $L$ satisfying
(\ref{dentwo}) over the set of pairs of circles $\big( J_1, J_2 \big)$
satisfying (\ref{denone}) is achieved when $J_1$ and $J_2$ touch at a
point, with $L$ being the line segment tangent to $J_1$ at the diametrically
opposed point and having endpoints in $J_2$ (see Figure $6$). For this
choice of line segment,
$\big\vert L \big\vert$ satisfies
\begin{eqnarray}
 \big\vert L \big\vert  & = & 2 \sqrt{ \rout \big( K \big)^2 - \Big( 2
 \rin \big( K \big) - \rout \big( K \big) \Big)^2 } \label{lenell} \\
 & = & 4  \sqrt{\rin \big( K \big) \Big(  \rout ( K )  -  \rin (
 K )  \Big)}, \nonumber
\end{eqnarray}
as required. $\Box$ \\
\noindent{\bf Proof of Corollary \ref{nthmone}}
From the first part of Lemma \ref{lemfi} follows the inclusion
\begin{eqnarray}
 & & \Big\{  \mathcal{L} \big( \conv B \big) > T^{\frac{5}{6} +
 \epsilon} \Big\}  \label{qtl} \\
 & \subseteq & \Big\{ \rin \big( B \big)
 > C_5 T \Big\} \cup \Big\{  \rout ( B )  -  \rin ( B ) \geq \frac{1}{16 C_5}
 T^{\frac{2}{3} + 2 \epsilon} \Big\}, \nonumber
\end{eqnarray}
for any fixed $C_5 > 0$.
Note that the event $Q_{f,g}^c$ occurs when the choices $f(T)=T$ and
$g(T)= \rin \big( B \big)$ are made. From
Lemma \ref{remaftdef}, it follows that
\begin{equation}\label{qtm}
 \mathbb{P} \Big( \rin \big( B \big) > C_5 T \Big) \leq  \exp \Big\{ - \frac{C_5^2}{128}
 T \Big\}.
\end{equation}
We deduce that
\begin{eqnarray}
 & &  \mathbb{P} \Big( \Big\{  \mathcal{L}\big( \conv B \big) > T^{\frac{5}{6} + \epsilon} \Big\} \cap
 \Big\{ \big\vert \enc(B) \big\vert \geq \pi T^2 \Big\} \Big)
 \nonumber \\
 & \leq &   \mathbb{P} \Big( \Big\{  \rout ( B )  -  \rin ( B ) > \frac{1}{16C_5}T^{\frac{2}{3} + 2\epsilon} \Big\} \cap
 \Big\{ \big\vert \enc(B) \big\vert \geq \pi T^2 \Big\} \Big) \nonumber
 \\
  &  & \qquad + \,  \mathbb{P} \Big( \rin \big( B \big) > C_5 T \Big) \nonumber
 \\
  & \leq & \exp \Big\{ - 2 \pi^2 T - c_5 T^{\frac{1}{3} + 4 \epsilon}
  \Big\} \, + \,  C \exp \Big\{ - \frac{C_5^2}{128} T \Big\}, \label{tty}
\end{eqnarray}
where $c_5$ is a constant satisfying $c_5 \leq c_4/(16 C_5)^2 $, and
where $\epsilon \in (0,1/{12})$. (Recall that the constant $c_4$
appeared in (\ref{xgh}).)
The first inequality is a consequence of (\ref{qtl}). In the
second, the bounds (\ref{qtm}) and (\ref{xgh}) were used. 
The choice $C_5 > 16 \pi$ ensures that the second
term in (\ref{tty}) is negligible. Thus, $c_5$ may be chosen to
be any value satisfying $c_5 \leq 1/\big(32(16)^4\big) = 2^{-21}$. 
The first part of the corollary
follows by applying Lemma \ref{ltwo}.

For the second part of the corollary, note that, if $\rout(B) \leq T
+ T^{\frac{2}{3} + \epsilon}$, then $\arcl \big( \partial( \conv B)
\big) \leq 2\pi \big( T + T^{\frac{2}{3} + \epsilon} \big)$, because
the arclength of the boundary of a planar convex set is monotone under
containment \cite{Santalo}. Thus, the second part of Theorem \ref{tgg} yields
the result. $\Box$\\

\noindent{\bf Acknowledgments:} We are grateful to Senya
Shlosman for proposing the model to us, and thank him in particular for
presenting heuristic arguments about the nature of its deviation.
We thank Manjunath Krishnapur, Ben Hough and Dapeng Zhan 
for helpful comments on a draft version of the paper.
\end{section}
\bibliography{loopbib}
\bibliographystyle{plain}
\end{document}